\documentclass{article}
\usepackage{graphicx} 
\usepackage{amsfonts}
\usepackage{amsmath}
\usepackage[numbers]{natbib}
\usepackage[letterpaper, margin=1in]{geometry}
\usepackage{cancel}
\usepackage{xcolor}
\usepackage{comment}
\usepackage{mathtools}
\usepackage{dsfont}
\usepackage{natbib}
\usepackage{amssymb}
\usepackage{enumitem}
\usepackage{amsthm}

\usepackage{color}
\usepackage{hyperref}

\title{Non-Ergodic Markovian Mirror Descent with Unbounded Noise}
\author{
\begin{minipage}{0.45\textwidth}
    \centering
    Yunus Alp Biyikoglu \\
    ECE and CSL \\
    UIUC \\
    \texttt{yunusab2@illinois.edu}
\end{minipage}
\hfill
\begin{minipage}{0.45\textwidth}
    \centering
    R. Srikant \\
    ECE, CSL, and NCSA \\
    UIUC \\
    \texttt{rsrikant@illinois.edu}
\end{minipage}
}

\begin{document}

\maketitle

\begin{abstract}
    We analyze stochastic mirror descent for minimizing a convex objective expressed as an expectation under the stationary distribution of a Markov chain, using noisy subgradient observations along a single sample path. Through Poisson’s equation, we decompose the Markovian noise and control its contribution under regularity and moment conditions that allow the stochastic subgradients and the Poisson solution to be unbounded across the state space. We establish an $\mathcal{O}(T^{-1/2})$ bound on the expected suboptimality gap of the weighted-average iterate for an appropriately chosen step size. The analysis does not require aperiodicity or convergence of the chain’s marginal distributions to stationarity. For irreducible finite-state chains, we verify the required conditions on the Poisson solution and derive explicit bounds that quantify the effect of periodicity. For irreducible countable-state chains, we verify the required regularity and moment conditions under a geometric Lyapunov drift condition and suitable growth bounds on the stochastic subgradients.
\end{abstract}

\section{Introduction}
Mirror descent, first introduced in \cite{nemirovskij1983problem}, is an algorithm employed to optimize  
an objective over a specified parameter set, akin to gradient or subgradient-based optimization methods. Since its inception, it has spurred several variants, most notably a nonlinear projected subgradient algorithm where a particularly chosen \emph{prox-function}, which induces a distance-like function denoted as \emph{Bregman divergence}, can influence the rate of convergence, especially in dimension-dependent constants \cite{beck2003mirror}. Utilizing an appropriate choice of a prox-function, it is possible to recover the iterates and convergence rate of what would otherwise be generated by a gradient or subgradient-based descent algorithm. However, other choices for the prox-function may yield a reduction of the upper bound on the suboptimality gap in dimension-dependent terms, such as when optimizing over a parameter set featuring non-Euclidean geometry, while preserving the $\mathcal{O}(T^{-1/2})$ rate found in (sub)gradient-type methods for general convex functions.

In certain problems, the observations of the objective itself or its (sub)gradient may be corrupted by noise conforming to a specified stochastic process. However, the algorithm may still guarantee convergence to the minimum value of the objective function, in expectation or with high probability. Such algorithms have been widely examined in the domain of gradient-type methods (see \cite{nemirovski2009robust, bottou2018optimization, nedic2001incremental, polyak1992acceleration} and the references therein). In the domain of stochastic mirror descent, several approaches to modeling randomness have been considered, such as independently and identically distributed (i.i.d.) noise sequences \cite{nemirovski2009robust,lan2012validation,lu2019relative, zhou2017stochastic, bayandina2018mirror}, as well as noise with zero bias \cite{hanzely2021fastest,nedic2014stochastic}. In general, noise sequences may not be unbiased at finite times, but may only be unbiased in an asymptotic limit. In a notable work \cite{duchi2012ergodic}, the authors consider \emph{ergodic} processes whose distributions converge to a stationary distribution over time and which admit a mixing time to characterize the rate of convergence. They show that their algorithm yields a convergence rate of $\mathcal{O}\left(\sqrt{\tau_{\text{mix}}/T}\right)$, where $\tau_{\text{mix}}$ is the mixing time of the random process. In the special case of irreducible, finite-state, and aperiodic Markov chains, this translates to a rate of $\mathcal{O}\left(\sqrt{(\log T)/T}\right)$, with coefficients that depend, among other things, on the mixing rate of the Markov chain. Two more recent works in which the mirror descent algorithm has been examined subject to a uniformly ergodic Markov chain are \cite{zhao2023markov,solodkin2408methods}, where in \cite{zhao2023markov} the author considers a finite-state, irreducible, and aperiodic chain whose stationary distribution is the $n$-point uniform distribution, and in \cite{solodkin2408methods} the authors allow for more general state spaces and characterize a mixing time to control noise terms. Another recent work on stochastic mirror descent under a process whose realizations across time are dependent is \cite{fort2026stochastic}, where the authors leverage a ``forgetting'' condition that is compatible with certain types of geometrically or polynomially decaying dependence, allowing certain classes of ergodic Markov chains. These works leave open the question of whether stochastic mirror descent converges for non-ergodic processes. We answer this question affirmatively for a large class of Markov chains.

A common approach to studying problems involving Markov chains is to use Poisson's equation; see \cite{haque2024stochastic, srikant2025rates, paul2026stochastic} for some recent works. Poisson's equation has been used to decompose Markovian noise terms into the sum of several terms, including a martingale difference sequence (MDS), which facilitates convergence proofs. Poisson's equation is an equation of the form
\begin{equation}
    u(x)-\mathbb{E}_{X\sim\pi}[u(X)]=w(x)-\mathbb{E}[w(X_1)|X_0=x],
\end{equation}
where $\{X_t \}_{t\in\mathbb{N}}$ is a Markov chain with state space $\mathcal{X}$ and stationary distribution $\pi$, and $u,w:\mathcal{X}\rightarrow\mathbb{R}^{n}$. Here, $u$ is assumed to be known and $w$ is said to be the solution of the equation. It is well known that the solution to Poisson's equation may exist even for non-ergodic Markov chains, which we exploit in this paper. However, to the best of our knowledge, except for specific examples, a characterization of the properties of the solution to the Poisson equation in terms of the parameters of the Markov chain in the context of convergence analysis of optimization algorithms does not seem to exist in the literature. One of the key contributions of our paper is such a characterization. 

A recent paper \cite{haque2024stochastic} studies the convergence of stochastic approximation algorithms with Markovian noise. They present their results in terms of the properties of the solution of Poisson equations and recognize that their results may hold even when the Markov chain is periodic (hence, non-ergodic). However, these properties of Poisson's equation are assumed (as in much of the prior literature) and are only verified for very specific examples of Markov chains. Here, we provide a framework to verify those assumptions and we also prove additional properties needed for the Markovian mirror descent algorithm.

In a nearly contemporaneous work \cite{paul2026stochastic}, the authors study mirror descent under Markovian noise. They prove asymptotic convergence of the algorithm for both convex and non-convex objectives and obtain a finite-time concentration bound for convex objectives. We believe that our work and their work provide complementary results. They allow for Markov chains that are parameterized by the iterate generated by mirror descent, which we do not consider. However, they assume that stochastic subgradients and the corresponding solutions of the Poisson equation are bounded over the state space of the chain, an assumption that need not be satisfied by countable state-space chains considered in our paper. Additionally, as in much of the prior work, their results are in terms of properties assumed for the solution of Poisson's equation. Here we prove all the properties that we require for irreducible finite-state Markov chains and for irreducible countable state-space Markov chains that satisfy a Lyapunov drift condition. Since our primary motivation is to show convergence to optimality, we do not consider non-convex objectives whereas \cite{paul2026stochastic} consider non-convex functions and establish convergence to some stationary point.

Our contributions in this paper are as follows:

\textbf{Performance of the Stochastic Mirror Descent Algorithm Under Markovian Noise:} We show that the objective function converges to its optimal value in expectation at a rate of $\mathcal{O}\left(T^{-1/2}\right)$, which represents an improvement over the convergence rate in \cite{duchi2012ergodic} and matches similar rates in \cite{fort2026stochastic} and in standard results in convex optimization where the objective is not assumed to be strongly convex. It should be noted that our results hold for certain classes of non-ergodic Markov chains while the results in \cite{duchi2012ergodic} require ergodicity and mixing, but allow more general processes than Markov chains. Additionally, unlike \cite{fort2026stochastic}, we do not require a forgetting condition or explicit decay of dependence over time, thereby accommodating certain classes of non-ergodic Markov chains, including periodic chains.

\textbf{Use of the Poisson Equation for Mirror Descent}: The use of Poisson's equation to study stochastic approximation is well known and has a long history \cite{benveniste2012adaptive}. To the best of our knowledge, this paper and \cite{paul2026stochastic} are the first to apply it to mirror descent. Mirror descent requires additional properties of Poisson's equation which we identify and prove in this paper.  For finite-state Markov chains, irreducibility is sufficient to ensure that a bounded solution to the Poisson equation exists. For countable state-space chains, we have utilized additional assumptions with respect to the drift behavior of the chain as well as the existence of upper bounds on the norm and Lipschitz coefficient of the gradient in terms of a specified Lyapunov function. We provide an example showing that our bound in terms of $\lVert Y\rVert_\infty$, where $Y$ is closely related to the fundamental matrix of a Markov chain, can in certain cases yield much-improved results over bounds derived using the mixing rate of the Markov chain, even for ergodic chains.

\textbf{An Explicit Derivation for Irreducible, Finite-State, and Periodic Markov Chains:} In Section 5.2, we decompose $\lVert Y\rVert_\infty$, a constant depending on the transition matrix of the Markov chain at hand, into its counterparts $\lVert Y_r\rVert_{\infty}, \, r\in \{0,1,...,d-1\},$ for the $d-$step transition process where $d$ denotes the period of the chain. To the best of our knowledge, this is a new result that explicitly incorporates the periodicity of the chain into the analysis to get a finite upper bound on the expected value of the suboptimality gap. Additionally, we show that this decomposition allows us to get mixing-rate-type bounds for periodic chains, where these mixing rates correspond to those of the Markov chain induced by the $d-$step transition process.

\textbf{Allowing Unbounded Subgradients:} In several previous works \cite{nemirovski2009robust, duchi2012ergodic}, the square of the norm of the (sub)gradient function is assumed to be bounded from above by a finite and state-independent constant in expectation, either unconditionally or conditionally with respect to a filtration. In \cite{lu2019relative}, the author generalizes this by employing a \emph{relative continuity} assumption, where the norm of the (sub)gradient is allowed to vary spatially with respect to the Bregman divergence $D_h$. In this work, we allow norms of subgradients and other related first or second-order quantities to be unbounded across the state space, but assume that they are upper bounded in expectation as a function of the initial state of the chain. These upper bounds can be made explicit in the case of finite-state chains and for countable state-space chains which admit a Lyapunov drift. This allows us to express the norm and Lipschitz coefficient of the subgradient and the solution to the Poisson equation as a function of the state of the Markov chain, which may be unbounded over the whole state space, unlike in \cite{paul2026stochastic} where these quantities are assumed to admit uniform upper bounds.

The paper is organized as follows: Section 2 begins with a discussion on a few relevant examples of mirror descent and then transitions to a particular application which is shown to admit a solution using the frameworks developed in the subsequent sections. Section 3 provides an explicit description of the problem at hand and lays out the assumptions we utilize in the subsequent sections. In Section 4, we present Theorem 4.1 and its corollary, which provides a suitable choice of the step size policy, to verify that a bound of  $\mathcal{O}\left(T^{-1/2}\right)$  on the expected value of the suboptimality gap can be achieved by our algorithm. In the following two sections, we focus on the cases of finite-state and countable state-space Markov chains, respectively, and verify the assumptions presented in Section 3 regarding the solution to the Poisson equation, with additional assumptions where needed. We include all proofs that are omitted in the main text in the Appendices.

\section{Motivating Applications in Machine Learning}

Reinforcement learning (RL) is a classical application of stochastic optimization in a Markovian framework. In RL, the goal is typically to find the optimal policy to be utilized by an agent such that it yields the maximum accrued rewards granted by the environment, along trajectories of Markov Decision Processes (MDPs). The nature of the dynamics of these processes warrants the creation of algorithms with provable convergence guarantees even when samples are subject to Markovian dependence. A recent work that makes strides in this area is \cite{haque2024stochastic}, where the authors develop a general-purpose framework for stochastic approximation with potentially unbounded Markovian noise and establish finite-time convergence rates, with applications to temporal-difference (TD) learning and Q-learning.

Methods based on stochastic mirror descent have also been explored in solving MDPs. \cite{jin2020efficiently, tiapkin2022primal} explore primal-dual stochastic mirror descent frameworks for solving several classes of MDPs, where they have access to unbiased observations of the gradient. Another important application of mirror descent in RL literature is in policy improvement, in which the non-Euclidean nature of probability simplices positions mirror descent as a favorable algorithm to update policies using estimates of the Q-function or the policy gradient while penalizing movements away from the previous policy through the Kullback-Leibler divergence, as touched upon in \cite{solodkin2408methods} and explored to a greater extent in the references therein. 

As we have previously discussed, we are interested in providing tools for general problems where observations of the objective or its subgradient are generated by a Markovian trajectory, similar to the problems discussed in \cite{solodkin2408methods, paul2026stochastic}. Since applications of mirror descent and its special case, gradient descent, to RL are well established, in the remainder of this section, we outline a graph learning problem that naturally admits a mirror descent formulation.

\subsection{Learning in Social Networks}
Here we show that the learning problem studied in \cite{katselis2018mixing, xie2021consistency} can be solved using our mirror descent framework. Consider a directed graph $\mathcal{G}=(\mathcal{V},\mathcal{E})$ with $|\mathcal{V}|=p$ nodes. Let $X_t\in \{0,1 \}^p$ denote the state vector of the process at time $t$, whose $i$th element $X_t(i)\in \{0,1 \}$ denotes the state of node $i$ at time $t$, for all $i \in \{1,2,...,p \}$ and all $t\in\mathbb{N}$. Conditioned on $X_t$ and another process, $W_t$, we assume that the elements of $X_{t+1}$ are generated independently of one another according to the following Bernoulli distribution:
\begin{equation}
    X_{t+1}(i)\sim\mathrm{Ber}(A(i,\cdot)X_{t} +b(i)W_{t}(i))
\end{equation}
where $A \in [0,1]^{p\times p}$, $b\in [b_{\text{min}},1]^p$ for some $b_{\text{min}}\in(0,1)$ (note that we do not require $b_{\text{min}}=\min_i b(i)$), and the elements of $W_{t}$ are generated independently of each other and of $X_t$ (and all other random variables belonging to times at or earlier than $t$) as $W_{t}(i)\sim \text{Ber}(1/2)$. We assume that $\sum_{j=1}^p A(i,j)+b(i)=1\,\,\,\forall\,i\in\{1,2,...,p\}$. We interpret the matrix $A$ to encode weights of directed edges between nodes on the graph. More explicitly, $A(i,j)$ denotes the weight of the directed edge from node $j$ to node $i$, with $A(i,j)=0\iff (j,i)\notin\mathcal{E}$. Among other things, the model can be used to represent opinion dynamics in graphs, where the state of a node represents a positive opinion ($1$) or a negative opinion ($0$) and the matrix $A$ encodes how different agents in the social network influence the evolution of each other's opinions.

It follows that the process $\{X_t \}_{t\in\mathbb{N}}$ corresponds to a finite-state, irreducible, and aperiodic Markov chain, whose state space $\mathcal{X}$ has a cardinality of $2^p$ and whose stationary distribution we denote by $\pi_X$. Our goal is to estimate the parameters $A$ and $b$ through the application of mirror descent on the negative log-likelihood function of the states $X_t$ along an observed trajectory of the chain. As in prior work, we assume that $b_{\text{min}}$ is known.

Let $\Omega\coloneqq\{\omega\in\ \mathbb{R}^{p\times (p+1)}\,|\, \omega(i,\cdot)\ge 0, \,\omega(i,p+1)\ge b_{\mathrm{min}}, \,\sum_{k=1}^{p+1}\omega(i,k)=1\,\,\forall\,i\in\{1,2,...,p\}\}$, which denotes the space of parameters we are interested in estimating by examining a sample Markovian trajectory $\{X_t \}$. Let $\Omega_r$ be the space in which the parameters belonging to a single node lie (corresponding to a row of a given element of $\Omega$), i.e. $\Omega_r\coloneqq\{\omega_r\in\ \mathbb{R}^{1\times (p+1)}\,|\, \omega_r\ge 0 
, \,\omega_r(p+1)\ge b_{\mathrm{min}}, \,\sum_{k=1}^{p+1}\omega_r(k)=1\}$.

For a given $\omega\in\Omega$, let $\hat{A}\in\mathbb{R}^{p\times p}$ and $\hat{b}\in \mathbb{R}^p$ denote sections of $\omega$ such that $\omega = [\hat{A}, \hat{b}]$. Additionally, let $\omega_i$ denote the $i$th row of $\omega$, where $\omega_i=[\hat{A}(i,\cdot),\hat{b}(i)]\in\Omega_r$.
Let $P_\omega(\cdot)$ denote the likelihood function of observations of the Markovian variables concerning all nodes, parameterized by $\omega$, and let $P^{(i)}_{\omega_i}(\cdot)$ be its analog for node $i$ of the graph, parameterized by $\omega_i$. For any two states $u$ and $v$ of the Markov chain and any $z\in\{0,1\}$, we have the following:
\begin{gather}
    P^{(i)}_{\omega_i}(z|u)=\left( \hat{A}(i,\cdot)u+\frac{\hat{b}(i)}{2}\right)^{z}\left( 1-\hat{A}(i,\cdot)u-\frac{\hat{b}(i)}{2}\right)^{1-z},\,\,\,\, P_\omega(v|u)=\prod_{i=1}^p P^{(i)}_{\omega_i}(v(i)|u).
\end{gather}
For an observed trajectory $\{X_0,X_1,...,X_{T-1},X_{T}\}$, we express the log-likelihood function $L_T(\omega)$ as follows:
\begin{align}
    L_T(\omega)=\frac{1}{T}\log P_\omega(X_0,X_1,...,X_{T-1},X_T)\nonumber
    &=\frac{1}{T}\left(\log P_\omega(X_0)+\sum_{t=1}^T \log P_\omega(X_t|X_{t-1})\right)\nonumber\\
    &=\frac{1}{T}\left(\log P_\omega(X_0)+\sum_{i=1}^p \sum_{t=1}^T\log P^{(i)}_{\omega_i}(X_t(i)|X_{t-1})\right)
\end{align}
We further assume that the distribution of $X_0$ is not a function of $\omega$, i.e. $P_\omega(X_0)\equiv P(X_0)$. Let $\omega_0\in\Omega$ be the actual but unknown parameters of the graph, i.e. $\omega_0=[A,b]$, and let $\omega^{\mathrm{ML}}_T\in\arg\min_{\omega\in\Omega} -L_T(\omega)$. Under the assumptions made up to this point, $\omega^{\mathrm{ML}}_T \xrightarrow[T\to\infty]{\text{a.s.}} \omega_0$, as proven in \cite{xie2021consistency}.

Notice that, ignoring the distribution of the initial state, the log-likelihood function is a sum of functions depending only on $\omega_i\in\Omega_r$ for each node $i$. This means that we can decouple the objective and optimize for the parameters associated with each node independently. Define a Markov chain $\{Y^{(i)}_t\}_{t\in\mathbb{N}}$ where $Y^{(i)}_t\coloneqq [X_{t};\,X_{t+1}(i)]$. It follows that $\mathcal{Y}=\mathcal{X}\times \{0,1\}$, where $\mathcal{Y}$ is the shared state space of the family of $\{Y^{(i)}_t \}$ chains indexed by $i$, since in the original $\{X_t \}$ chain a state with a given value among $\{0,1\}$ for its $i$th element can be reached from any other state in one step with positive probability. Since the original chain has a positive one-step probability between any two states, the chain $\{Y^{(i)}_t \}$ is irreducible. Our objective of the maximum-likelihood (ML) problem associated with node $i$ becomes $-\frac{1}{T}\sum_{t=0}^{T-1}\log P_{\omega_i}^{(i)}(Y^{(i)}_{t}(p+1)\,|\,Y^{(i)}_{t}(1:p))$. Note that this quantity is convex in $\omega_i$. As $T\rightarrow\infty$, the objective of the ML problem converges pointwise in $\omega_i$ to $\mathbb{E}_{Y\sim\pi_{Y^{(i)}}}[-\log P_{\omega_i}^{(i)}(Y^{(i)}(p+1)\,|\,Y^{(i)}(1:p))]$ almost surely, by the ergodic theorem for Markov chains relating to time averages (which does not require aperiodicity), where $\pi_{Y^{(i)}}(y)\coloneqq \pi_X(y(1:p))P^{(i)}_{(\omega_0)_i}(y(p+1)\,|\,y(1:p))$ is the stationary distribution of $\{Y^{(i)}_t\}$. It also follows that $(\omega_0)_i$ is the unique minimizer of $\mathbb{E}_{Y^{(i)}\sim\pi_{Y^{(i)}}}[-\log P_{\omega_i}^{(i)}(Y^{(i)}(p+1)\,|\,Y^{(i)}(1:p))]$, as we show in Appendix A.

We will define an invertible transformation to account for the fact that every $\omega_i\in\Omega_r$ has an immovable mass $b_{\text{min}}$ in its last coordinate. Let $\Delta=\{u\in\mathbb{R}^{p+1}\,|\,u_i\ge 0\,\,\forall\, 1\le i\le p+1, \mathds{1}^\top u=1 \}$ denote the probability simplex in $\mathbb{R}^{p+1}$ and let the transformation $\Phi:\Omega_r\rightarrow\Delta$ be as follows:
\begin{equation}
    \Phi(\omega_i)=\frac{1}{1-b_{\text{min}}}[\omega_i(1:p)^\top;\,\omega_i(p+1)-b_{\text{min}}]\,\,\,\forall\,\omega_i\in\Omega_r.
\end{equation}
It then follows that the inverse of $\Phi$ is given by
\begin{equation}
    \Phi^{-1}(u)=(1-b_{\text{min}})[u(1:p)^\top,\,u(p+1)+b_{\text{min}}/(1-b_{\text{min}})]\,\,\,\forall\,u\in\Delta.
\end{equation}

Let $f^{(i)}(\omega_i, Y^{(i)})\coloneqq -\log P^{(i)}_{\omega_i}(Y^{(i)}(p+1)\,|\,Y^{(i)}(1:p))$. Noting that there is a one-to-one mapping between $\Omega_r$ and $\Delta$, we define the objective $\bar{f}_M^{(i)}$ as follows:
\begin{gather}
    f_M^{(i)}(u,y)\coloneqq f^{(i)}(\Phi^{-1}(u),y) \,\,\forall\,u\in\Delta, y\in\mathcal{Y},\nonumber\\
    \bar{f}_M^{(i)}(u)\coloneqq\mathbb{E}_{Y^{(i)}\sim\pi_{Y^{(i)}}}[f_M^{(i)}(u,Y^{(i)})]\,\,\forall\,u\in\Delta.
\end{gather}

Hence, our goal is to minimize $\bar{f}^{(i)}_M(u)$ where $u \in\Delta$, using observations of the stochastic gradient $g^{(i)}_M(u,y)\coloneqq \nabla_u f^{(i)}_M(u,y)\,\,\,\forall \,u\in \operatorname{relint}(\Delta), \,y\in\mathcal{Y}$, whose stationary expectation under $\pi_{Y^{(i)}}$ is a gradient of the objective $\bar{f}^{(i)}_M$. We define a mirror descent iterate sequence $\{u^{(i)}_t\}_{t\in\mathbb{N}}$ where $u^{(i)}_0 = (1/(p+1))\mathds{1}$ and update $u^{(i)}_{t+1}$ from $u^{(i)}_t$ as follows, for a suitably chosen step size $\epsilon_t$:
\begin{gather}
    u^{(i)}_{t+1}= \arg\min_{u \in \Delta}\, \langle g_M^{(i)}(u^{(i)}_t, Y^{(i)}_t), u \rangle + \frac{1}{\epsilon_t}D_{\text{KL}}(u\| u^{(i)}_{t})\nonumber\\
    \implies u^{(i)}_{t+1}(j)=\dfrac{u^{(i)}_{t}(j)\exp(-\epsilon_t\,g_M^{(i)}(u^{(i)}_t,Y^{(i)}_t)(j))}{\sum_{k=1}^{p+1} u^{(i)}_{t}(k)\exp(-\epsilon_t\,g_M^{(i)}(u^{(i)}_t,Y^{(i)}_t)(k))}\,\,\,\forall\, 1\le j\le p+1,
\end{gather}
where $D_{\text{KL}}(u\|v)\coloneqq \sum_{j=1}^{p+1}u_j\log(u_j/v_j)$ denotes the Kullback-Leibler divergence between probability vectors $u$ and $v$.

It follows that this problem satisfies Assumptions 3.1-3.5. Invoking Theorem 4.1 and utilizing several other arguments, we obtain the following result: $\bar{u}^{(i)}_T\coloneqq(\sum_{t=0}^T \epsilon_t u^{(i)}_t)/(\sum_{t=0}^T \epsilon_t)\xrightarrow[T\to\infty]{p}(u^{(i)})^\star\coloneqq \Phi((\omega_0)_i)$ 
for all $i \in \{1,2,...,p\}$. We provide more details in Appendix A.

\section{Problem Description and Assumptions}
In this paper, we analyze the Markovian mirror descent (MMD) algorithm for solving a class of nonsmooth convex optimization problems in which the function being optimized is formulated as the stationary expectation of a series of functions that are corrupted by Markovian noise. More formally, let $\{f(\cdot,x), \, x \in \mathcal{X} \}$ be a collection of functions over a domain that includes the convex and compact set $\Theta \subseteq \mathbb{R}^n$, where $\mathcal{X}$ denotes the state space of a discrete-time, time-homogeneous Markov chain with unknown stationary distribution $\pi(\cdot)$. We are interested in minimizing the convex function $\bar{f}(\theta) \coloneqq \mathbb{E}_{X \sim \pi}[f(\theta,X)]$ over $\Theta$. We assume there exist vector-valued functions $\{g(\cdot,x), \, x \in \mathcal{X} \}$ whose stationary expectation, $\bar{g}(\theta)\coloneqq\mathbb{E}_{X \sim \pi}[g(\theta,X)]$, is a subgradient of $\bar{f}$ evaluated at the point $\theta$, i.e. $\bar{g}(\theta) \in \partial\bar{f}(\theta)$, for all $\theta$ belonging to a set in which the iterates of the algorithm we describe below lie, in a manner we will make precise later.

Assume that there exists a discrete-time Markov chain $\{X_t\}_{t \in \mathbb{N}}$ from which the stochastic subgradients $g(\cdot,X_t)$ are generated at each time step $t$. Beginning with a predetermined initial parameter value, $\theta_0$, we propose the following MMD algorithm to update $\theta_{t+1}$ for $\,t \in \mathbb{N}$, where $\langle \cdot\,, \cdot\rangle$ denotes the standard Euclidean inner product and $\{\epsilon_t\}_{t \in \mathbb{N}} $ denotes the positive and non-increasing step size sequence:
\begin{equation}
    \label{eq1}
    \theta_{t+1} = \arg\min_{\theta \in \Theta}\, \langle g(\theta_{t}, X_{t}),\theta \rangle + \frac{1}{\epsilon_t}D_h(\theta, \theta_{t})
\end{equation}
Here, $D_h(\theta_1,\theta_2)$ denotes the Bregman divergence between $\theta_1$ and $\theta_2$ and is defined in terms of the $\mu$-strongly convex function $h$ as follows (where $\mu > 0$):
\begin{equation}
\label{dh}
    D_h(\theta_1,\theta_2) \coloneqq h(\theta_1)-h(\theta_2)-\langle \nabla h(\theta_2), \,\theta_1-\theta_2 \rangle 
    \ge \dfrac{\mu}{2}\lVert \theta_1-\theta_2\rVert^2
\end{equation}
\textbf{Assumption 3.1.} 
$\{X_t\}_{t \in \mathbb{N}}$ is a discrete-time, time-homogeneous Markov chain with state space $\mathcal{X}$ and unique stationary distribution $\pi$.

\textbf{Assumption 3.2.} The functions that define the problem at hand and their domains satisfy the following:
\begin{enumerate}[label=(\roman*)]
    \item $\Theta \subset \mathbb{R}^n$ is a convex and compact set. Let $\lVert \cdot\rVert$ be any norm in $\mathbb{R}^n$ and let $\lVert \cdot \rVert_*$ denote its associated dual norm, where $\lVert y \rVert_* = \sup_{\lVert z\rVert\le 1} \langle y,z\rangle$. Since $\Theta$ is compact, there exist constants $R, R_*>0$ such that $\lVert \theta_2 - \theta_1 \rVert \le R,\,\, \lVert \theta_2 - \theta_1 \rVert_* \le R_*\, \, \forall \, \theta_1, \theta_2 \in \Theta$.
    \item There exist functions $f:\Theta_f \times \mathcal{X}\rightarrow \mathbb{R}$,\, $g: \Theta_g\times \mathcal{X} \rightarrow \mathbb{R}^n$, and $h: \Theta_h \rightarrow \mathbb{R}$ where $\Theta \subseteq \Theta_f \cap \Theta_h\,$ and $\,\Theta_f \cup \Theta_g \cup \Theta_h \subseteq \mathbb{R}^n$. 
    \item $h$ is 
    differentiable on a domain containing $\operatorname{relint}(\Theta)$ and is $\mu$-strongly convex over $\Theta$ for some $\mu > 0$. 
    \item There exists a set $\widetilde{\Theta}$ such that $\operatorname{relint}(\Theta) \subseteq \widetilde{\Theta} \subseteq\Theta\,\cap\,\Theta_g\,\cap\,\text{dom}(\nabla h).$ This set satisfies several other properties as outlined below.
    \item The inequality (\ref{dh}) holds for all $\theta_1 \in \Theta$ and all $\theta_2 \in \widetilde{\Theta}$.
    \item Let $\,\bar{f}(\theta)\coloneqq \mathbb{E}_{X\sim \pi}[f(\theta,X)] \,\,\, \, \forall\,\theta \in \Theta_f$, $\,\,\,\bar{g}(\theta) \coloneqq \mathbb{E}_{X\sim\pi}[g(\theta,X)]\,\,\,\,\forall \, \theta \in \Theta_g,\,$ with \\$\,\,\mathbb{E}_{X\sim\pi}\left[\, \lVert g(\theta,X)\rVert_* \right] < \infty\,\,\,\forall \, \theta \in \Theta_g$. We assume $\bar{f}$ is finite-valued and convex over $\Theta$ and that for all $\theta \in \widetilde{\Theta}$, $\bar{g}(\theta)\in \partial\bar{f}(\theta)$. This implies the following inequality for all $\theta_1\in \widetilde{\Theta}$ and all $\theta_2 \in \Theta$:
    \begin{equation}
\label{eq3}
    \bar{f}(\theta_2) - \bar{f}(\theta_1) \ge  \langle \bar{g}(\theta_1),\,\theta_2 - \theta_1 \rangle.
\end{equation}
    \item The set of minimizers of $\bar{f}$ over $\Theta$ is nonempty. Let $\theta^\star$ denote such a minimizer.
    \item The function $h$ is such that, given any $\theta_t \in \widetilde{\Theta}, \,t\in \mathbb{N},$ and the parameter update rule in (\ref{eq1}), the resulting $\theta_{t+1}$ exists, is unique, and is also an element of $\widetilde{\Theta}$. A sufficient condition on $h$ to satisfy this assumption where $\widetilde{\Theta}=\operatorname{relint}(\Theta)$ is given in \cite[Equation 3.12]{beck2003mirror}, though we do not impose that condition here.
    \item $\theta_0 \in \widetilde{\Theta}$ is chosen to satisfy $D_h(\theta,\theta_0) \le M$ for all $\theta \in \Theta$, for some $M \in (0,\infty)$. If there exists a $\theta'\in \arg \min_{\theta \in \widetilde{\Theta}} h(\theta)$ where $\nabla h(\theta')=0$, then
    one possible way to ensure this is to set $\theta_0=\theta'$ with the corresponding choice for $M$ being $\sup_{\theta\in \Theta} h(\theta)-\ h(\theta')$ when this supremum is finite, similar to how it is taken to be in \cite{nemirovski2009robust}.
    \item For a sample path $\{X_t \}_{t \in \mathbb{N}}$ of the Markov chain at hand, the stochastic subgradient vectors $g(\cdot,X_t)$ are made available to the user at each time step $t$. Together with a suitably chosen $\theta_0$ and step size sequence $\{\epsilon_t\}_{t\in\mathbb{N}}$, the iterates $\{\theta_t \}_{t \in \mathbb{N}}$ are calculated using the update equation (\ref{eq1}).
\end{enumerate}

\emph{Remark}. Note that we have not assumed the functions $f(\theta,x)$ where $ x\in \mathcal{X}$ to be convex in $\theta$, nor have we assumed an explicit relation between $f(\theta,x)$ and $g(\theta,x)$ for any $x\in \mathcal{X}$, unlike in \cite{duchi2012ergodic}. However, if for a specified problem the functions $f(\theta,x)$ happen to be convex for all $x \in \mathcal{X}$, then the stochastic subgradient function $g$ may be chosen to satisfy $g(\theta,x) \in \partial_\theta f(\theta,x)\,\, \forall \, \theta \in \widetilde{\Theta}, \forall\, x \in \mathcal{X}$, assuming the swapping of the subgradient and expectation operations is valid.

\textbf{Assumption 3.3.} There exists a function $A_1: \mathcal{X} \rightarrow [0,\infty)$ such that the stochastic subgradient $g$ satisfies the following for all $\theta \in \widetilde{\Theta}$ and $x \in \mathcal{X}$:
\begin{equation}
\lVert g(\theta, x)\rVert_* \le A_1(x). 
\end{equation}

Moreover, there exist functions $r:\Theta_g\rightarrow \mathbb{R}^n$ and $A_2: \mathcal{X} \rightarrow [0,\infty)$ such that $\tilde{g}(\theta,x)\coloneqq g(\theta,x)-r(\theta)$ satisfies the following for all $\theta, \theta' \in \widetilde{\Theta}$ and $x \in \mathcal{X}$:
\begin{equation}
    \lVert \tilde{g}(\theta', x)-\tilde{g}(\theta,x) \rVert \le A_2(x)\lVert\theta' - \theta\rVert.
\end{equation}

\emph{Remark.} Through this assumption, we allow $g$ to be non-Lipschitz in $\theta$ through its additive component $r(\theta)$. When applicable, a suitable choice for $r(\theta)$ may be $\bar{g}(\theta)$. However, since $\bar{g}(\theta)$ might not be known or available in closed form, we present this more general condition which may be easier to verify in practice. 

\textbf{Assumption 3.4.} Let $\mathbb{P}_x(\cdot)$ denote the probability law of the chain initialized at $X_0=x$, and let $\mathbb{E}_x[\cdot]$ be the expectation under $\mathbb{P}_x$, for all $x\in\mathcal{X}$. Then, we assume there exists a function $V_\theta(x)$ for all $\theta \in \widetilde{\Theta}$ and $x \in \mathcal{X}$ that satisfies the Poisson equation:
\begin{equation}
    \label{eq4}
    g(\theta, x) - \bar{g}(\theta) = V_\theta(x) - \mathbb{E}_x\left[V_\theta(X_1)\right]
\end{equation}

Additionally, we assume that there exist functions $B_1, B_2:\mathcal{X}\rightarrow[0,\infty)$ such that the solution $V$ satisfies the following for all $\theta, \theta'\in \widetilde{\Theta}$ and $x \in \mathcal{X}$:

\begin{gather}
    (i)\,\,\lVert V_{\theta}(x)\rVert_* \le B_1(x), \nonumber\\
        (ii)\,\,\lVert V_{\theta'}(x)-V_{\theta}(x) \rVert \le B_2(x)\lVert\theta' - \theta \rVert. \, \,
\end{gather}

\emph{Remark.} We allow for the functions $A_1, A_2, B_1, B_2$ to be unbounded over the state space $\mathcal{X}$, as long as they are consistent with our other assumptions. This is a relaxation of the corresponding assumptions of \cite{paul2026stochastic}, which stipulate that there exist uniform bounds across both the parameter and state spaces on the norm and Lipschitz coefficient of the stochastic subgradient and solution to the Poisson equation. However, if the state space of the Markov chain is finite, then it follows that such uniform upper bounds exist in our case as well.

\textbf{Assumption 3.5.}
Let $\{X_t \}_{t\in \mathbb{N}}$ be a sample path of the Markov chain, with an initial state of $X_0=x_0$ for some $x_0\in \mathcal{X}$
. Then, we assume there exist functions $C_1,C_2,C_3,C_4:\mathcal{X}\rightarrow [0,\infty)$ such that, for all $t \in \mathbb{N}$, the following are satisfied:
\begin{gather}
    (i)\,\,\mathbb{E}_{x_0}[B_1(X_t)]\le C_1(x_0),\nonumber\\
    (ii)\,\,\mathbb{E}_{x_0}[A_1(X_t)^2] \le C_2(x_0), \nonumber \\
    (iii)\,\,\mathbb{E}_{x_0}[A_1(X_t)B_1(X_{t+1})] \le C_3(x_0), \nonumber \\
    (iv)\,\,\mathbb{E}_{x_0}[A_1(X_t)B_2(X_{t+1})] \le C_4(x_0). 
\end{gather}

\emph{Remark.} As we will verify in Section 5, a solution to the Poisson equation always exists for irreducible Markov chains with a finite state space, and the relevant quantities admit upper bounds over the finite set $\mathcal{X}$. In Section 6, by utilizing Assumption 3.3 and some additional properties, such as the existence of a Lyapunov function that satisfies a certain drift condition, we will show the existence of a valid solution in the case of countable state space chains and will derive expressions for $B_1(x),B_2(x),C_1(x),C_2(x),C_3(x),$ and $C_4(x)$.

\section{Convergence Analysis}

In this section, we present our main theorem on the convergence of the MMD algorithm in expectation, under the conditions described in the earlier sections.

\textbf{Theorem 4.1.} \emph{Consider the MMD algorithm presented in (\ref{eq1}). Let $\bar{\theta}_T \coloneqq \dfrac{\sum_{t=0}^T \epsilon_t \theta_t}{\sum_{t=0}^T \epsilon_t}$. With the assumptions laid out in Section 3, the suboptimality gap between $\bar{f}(\bar{\theta}_T)$ and $\bar{f}(\theta^\star)$ satisfies the following bound in expectation:}
\begin{equation}
    \label{convergence}
    0\le \mathbb{E}_{x_0}[\bar{f}(\bar{\theta}_T)] - \bar{f}(\theta^\star) \le \dfrac{1}{\sum_{t=0}^T \epsilon_t}\Biggl[M + 2C_1(x_0)R\,\epsilon_0 + \frac{1}{\mu}\left(\frac{C_2(x_0)}{2}+C_3(x_0)+R_*C_4(x_0)\right)\sum_{t=0}^{T}\epsilon_t^2 \Biggr]
\end{equation}
\begin{proof}
From the parameter update rule presented in (\ref{eq1}), note that $\theta_{t+1}$ is the solution to a convex optimization problem in $\theta$ over a compact set. Hence, the following optimality condition holds:

\begin{equation}
\label{eq6}
    \langle g(\theta_t, X_t) - \frac{1}{\epsilon_t}\nabla h(\theta_t) + \frac{1}{\epsilon_t}\nabla h(\theta_{t+1}),\, \theta - \theta_{t+1} \rangle \ge 0 \, \, \, \forall \, \, \theta \in \Theta, \, \forall \, t\in \mathbb{N}
\end{equation}
Inserting $\theta = \theta_t$ into  (\ref{eq6}), combining it with  (\ref{dh}), and using the definition of the dual norm, we get the following:
\begin{align}
\label{onestepbound}
     \epsilon_t\lVert  g(\theta_t, X_t)\rVert_* \lVert \theta_{t+1}-\theta_{t}\rVert\ge \epsilon_t\langle g(\theta_t, X_t),\,\theta_t-\theta_{t+1}\rangle &\ge\langle \nabla h(\theta_t),\,\theta_t-\theta_{t+1}\rangle+\langle\nabla h(\theta_{t+1}),\,\theta_{t+1}-\theta_{t}\rangle\nonumber  \\
    &= D_h(\theta_{t+1}, \theta_{t})+D_h(\theta_{t},\theta_{t+1}) \ge \mu \lVert \theta_{t+1}-\theta_t\rVert^2 \nonumber \\
    \implies \lVert \theta_{t+1} - \theta_t\rVert \le \frac{\epsilon_t }{\mu}\lVert  g(\theta_t, X_t)\rVert_* &\le \frac{\epsilon_t A_1(X_t)}{\mu} \, \, \forall \, t \in \mathbb{N},
\end{align}
where the last inequality follows from Assumption 3.3.

We proceed with the analysis as follows:
\begin{align*}
D_h(\theta^\star, & \,\theta_{t+1}) - D_h(\theta^\star, \theta_t) \\ &= -h(\theta_{t+1})-\langle \nabla h(\theta_{t+1}),\,\theta^\star - \theta_{t+1}\rangle + h(\theta_t) + \langle \nabla h(\theta_t),\,\theta^\star - \theta_t \rangle        &&  \\
     &\le  -h(\theta_{t+1}) + h(\theta_t) + \langle \nabla h(\theta_t),\,\theta_{t+1}-\theta_t\rangle + \epsilon_t\langle g(\theta_t, X_t),\,\theta^\star - \theta_{t+1}\rangle      && \text{(By (\ref{eq6}))} 
     \\
     &\le -\frac{\mu}{2} \lVert \theta_{t+1} - \theta_t \rVert^2  + \epsilon_t\langle g(\theta_t, X_t),\,\theta^\star - \theta_{t+1}\rangle              && \text{(By (\ref{dh}))} 
     \\
     &= -\frac{\mu}{2} \lVert \theta_{t+1} - \theta_t \rVert^2  + \epsilon_t\langle g(\theta_t, X_t),\,\theta^\star - \theta_{t}\rangle +  \epsilon_t\langle g(\theta_t, X_t),\,\theta_t - \theta_{t+1}\rangle \\
     &\le -\frac{\mu}{2} \lVert \theta_{t+1} - \theta_t \rVert^2  + \epsilon_t\langle g(\theta_t, X_t),\,\theta^\star - \theta_{t}\rangle  + \epsilon_t \lVert g(\theta_t, X_t) \rVert_*\lVert \theta_t - \theta_{t+1}\rVert   \\
     & \leq -\frac{\mu}{2} \lVert \theta_{t+1} - \theta_t \rVert^2  + \epsilon_t\langle g(\theta_t, X_t),\,\theta^\star - \theta_{t}\rangle  + \frac{\epsilon_t^2}{2\mu}\lVert g(\theta_t, X_t)\rVert_*^2 + \frac{\mu}{2} \lVert \theta_{t+1} - \theta_t \rVert^2 && \text{$\left(\sqrt{a^2b^2} \leq \frac{a^2 + b^2}{2}\right)$} 
     \\
     &\le \epsilon_t\langle g(\theta_t, X_t),\,\theta^\star - \theta_{t}\rangle  + \frac{\epsilon_t^2}{2\mu}A_1^2(X_t) && \text{(Assumption 3.3)}
     \\
     &= \epsilon_t \langle \bar{g}(\theta_t),\,\theta^\star - \theta_t\rangle + \epsilon_t \langle g(\theta_t, X_t)-\bar{g}(\theta_t),\,\theta^\star - \theta_t \rangle + \frac{\epsilon_t^2}{2\mu}A_1^2(X_t) \\
     &\le \epsilon_t (\bar{f}(\theta^\star) - \bar{f}(\theta_t)) + \epsilon_t \langle g(\theta_t, X_t)-\bar{g}(\theta_t),\,\theta^\star-\theta_t \rangle + \frac{\epsilon_t^2}{2\mu}A_1^2(X_t) && \text{(Assumption 3.2)}
\end{align*}

Rearranging the inequality, summing from $t = 0$ to $T$, and noting that the $D_h(\theta^\star, \theta_{t+1}) - D_h(\theta^\star, \theta_t)$ term telescopes, we get the following:
\begin{align*}
    \sum_{t=0}^{T} \epsilon_t (\bar{f}(\theta_t) - \bar{f}(\theta^\star)) &\le D_h(\theta^\star, \theta_0) - D_h(\theta^\star,\theta_{T+1}) + \frac{1}{2\mu}\sum_{t=0}^{T} A_1^2(X_t)\epsilon_t^2 + \sum_{t=0}^{T}\epsilon_t \langle g(\theta_t, X_t)-\bar{g}(\theta_t),\,\theta^\star -\theta_t\rangle\\
    &\leq M + \frac{1}{2\mu}\sum_{t=0}^{T} A_1^2(X_t)\epsilon_t^2 + \sum_{t=0}^{T}\epsilon_t \langle g(\theta_t, X_t)-\bar{g}(\theta_t),\,\theta^\star -\theta_t\rangle 
\end{align*}
By the convexity of $\bar{f}$, $\bar{f}(\bar{\theta}_T) \le \dfrac{\sum_{t=0}^T \epsilon_t \bar{f}(\theta_t)}{\sum_{t=0}^T \epsilon_t}$. Hence, by dividing our main inequality by $\sum_{t=0}^T \epsilon_t$ and taking expectation, we get
\begin{equation}
\label{eq8}
    0 \le \mathbb{E}_{x_0}\left[\bar{f}(\bar{\theta}_T)\right] - \bar{f}(\theta^\star) \le \dfrac{1}{\sum_{t=0}^T \epsilon_t}\left[M + \frac{C_2(x_0)}{2\mu}\sum_{t=0}^{T}\epsilon_t^2 + \mathbb{E}_{x_0}\left[ \underbrace{\sum_{t=0}^{T}\epsilon_t \langle g(\theta_t, X_t)-\bar{g}(\theta_t),\,\theta^\star -\theta_t\rangle}_{(a)} \right] \right]
\end{equation}
All that remains is to bound $(a)$ in expectation. Let $\mathcal{F}_t = \sigma(X_0, X_1, ..., X_{t-1}, X_t)$. Invoking Assumption 3.4, we get the following in the $\mathbb{P}_{x_0}-$a.s. sense:
\begin{align*}
    (a)=\sum_{t=0}^T \epsilon_t \langle g(\theta_t, X_t)-\bar{g}(\theta_t), \,\theta^\star-\theta_t \rangle&= \sum_{t=0}^T\epsilon_t \langle V_{\theta_t}(X_t)-\mathbb{E}_{x_0}\left[ V_{\theta_t}(X_{t+1})| \mathcal{F}_t \right], \,\theta^\star-\theta_t\rangle\\
    &= \sum_{t=0}^T\epsilon_t\langle V_{\theta_t}(X_t) - V_{\theta_{t+1}}(X_{t+1}), \,\theta^\star-\theta_t \rangle&& (i) \\
    &+\sum_{t=0}^T\epsilon_t \langle V_{\theta_{t+1}}(X_{t+1}) - V_{\theta_{t}}(X_{t+1}), \,\theta^\star-\theta_t \rangle && (ii) \\
    &+\sum_{t=0}^T\epsilon_t \langle V_{\theta_{t}}(X_{t+1}) -  \mathbb{E}_{x_0}\left[ V_{\theta_t}(X_{t+1})|\mathcal{F}_t \right], \,\theta^\star-\theta_t\rangle && (iii)
\end{align*}

Note that in (\ref{eq4}), the expectation is taken over the distribution $\mathbb{P}_x(X_1\in\cdot\,)$ and $\theta$ is treated as a non-random quantity. Since we wish to insert $\theta = \theta_t$ into (\ref{eq4}), the correct way to denote the expectation is to condition it on $\mathcal{F}_t$, since $\theta_t$ is measurable with respect to $\mathcal{F}_{t-1}$ and, given $\mathcal{F}_t$, $X_{t+1}$ is independent of $\theta_t$ and its distribution is fully characterized by $X_t$.

We bound $(i), (ii),$ and $(iii)$, as follows:
\begin{align}
    (i)&=\sum_{t=0}^{T} \epsilon_t\langle V_{\theta_t}(X_t)-V_{\theta_{t+1}}(X_{t+1}), \,\theta^\star-\theta_t\rangle \nonumber\\
    &= \epsilon_0 \langle V_{\theta_0}(X_0),\,\theta^\star - \theta_0\rangle- \epsilon_{T} \langle V_{\theta_{T+1}}(X_{T+1}),\,\theta^\star - \theta_T\rangle + \sum_{t=1}^{T}\langle V_{\theta_t}(X_t),\,\epsilon_t(\theta^\star - \theta_t) - \epsilon_{t-1}(\theta^\star-\theta_{t-1})\rangle \nonumber\\
    &\le R(\epsilon_0 B_1(X_0) + \epsilon_{T}B_1(X_{T+1})) +  \sum_{t=1}^{T}\, B_1(X_t) \lVert \epsilon_t(\theta^\star - \theta_t) - \epsilon_{t-1}(\theta^\star-\theta_{t-1}) \rVert
    \nonumber\\
    &= R(\epsilon_0 B_1(X_0) + \epsilon_{T}B_1(X_{T+1})) + \sum_{t=1}^{T} \, B_1(X_t) \lVert \epsilon_t(\theta_{t-1} - \theta_{t})+\epsilon_t(\theta^\star-\theta_{t-1}) - \epsilon_{t-1}(\theta^\star-\theta_{t-1}) \rVert \nonumber\\
    &\le R(\epsilon_0 B_1(X_0) + \epsilon_{T}B_1(X_{T+1})) + \sum_{t=1}^{T} \, B_1(X_t)  \left[ \epsilon_t \lVert\theta_t - \theta_{t-1}\rVert + R(\epsilon_{t-1}-\epsilon_t)\right] \nonumber\\
    &\le R(\epsilon_0 B_1(X_0) + \epsilon_{T}B_1(X_{T+1})) + \sum_{t=0}^{T-1} \, B_1(X_{t+1})  \left[ \epsilon_{t}^2 \frac{A_1(X_t)}{\mu} + R(\epsilon_{t}-\epsilon_{t+1})\right] \,\,\,\,\,\,\,\,\,\,\,\,\,\,\,\text{(By (\ref{onestepbound}))}
\end{align}

Taking expectation and invoking Assumption 3.5, we get
\begin{equation}
    \label{i} 
    \mathbb{E}_{x_0}[(i)] \le 2C_1(x_0)R\epsilon_0+ \frac{C_3(x_0)}{\mu}\sum_{t=0}^T\epsilon_t^2
\end{equation}
\begin{gather}
   (ii) \le R_*  \sum_{t=0}^T\epsilon_t B_2(X_{t+1})\lVert \theta_{t+1} - \theta_t \rVert \le \frac{R_*}{\mu}  \sum_{t=0}^T \epsilon_t^2 A_1(X_{t})B_2(X_{t+1}) 
    \implies \mathbb{E}_{x_0}[(ii)]\le \frac{R_*C_4(x_0)}{\mu}  \sum_{t=0}^T \epsilon_t^2 \label{ii}
\end{gather}

For term $(iii)$, note that the sequence $M_{t}\coloneqq \langle V_{\theta_{t-1}}(X_{t}) -  \mathbb{E}_{x_0}\left[ V_{\theta_{t-1}}(X_{t})|\mathcal{F}_{t-1} \right], \,\theta^\star-\theta_{t-1}\rangle$ for $t
\in \mathbb{Z}^+$ is a martingale difference sequence with respect to the filtration $\mathcal{F}_t$, where $\mathbb{E}_{x_0}[M_{t}|\mathcal{F}_{t-1}]=0$ (the integrability condition follows from Assumptions 3.4 and 3.5, with respect to the measure $\mathbb{P}_{x_0}$). Hence, we get
\begin{align}
\label{iii E}
    \mathbb{E}_{x_0}[ (iii) ] =& \sum_{t=0}^T\epsilon_t \,\mathbb{E}_{x_0}[M_{t+1}]= \sum_{t=0}^T\epsilon_t \,\mathbb{E}_{x_0}\left[\,\mathbb{E}_{x_0}[M_{t+1}|\mathcal{F}_t ]\,\right] = 0.
\end{align}

From here, combining (\ref{eq8}), (\ref{i}), (\ref{ii}), and (\ref{iii E}) yields the claimed result.
    
\end{proof}
\textbf{Corollary 4.2.} \emph{We propose two suitable choices of the step size sequence $\{\epsilon_t \}_{t\in\mathbb{N}}$:
\begin{enumerate}[label=(\roman*)]
\item If the horizon $T$ is unknown ahead of time, then a step size policy of $\epsilon_t=\dfrac{\alpha}{\sqrt{t+1}} \,\,\forall\,t\in\mathbb{N}$ where $\alpha > 0$ is suitable, yielding a convergence rate of $\mathcal{O}\left(\dfrac{\log\,T}{\sqrt{T}}\right)$ in expectation.
\item If the horizon $T$ is known ahead of time, then a constant step size policy can be used with $\epsilon_t=\dfrac{\alpha}{\sqrt{T+1}}\,\,\forall\,t\in\mathbb{N}$ where $\alpha > 0$, yielding a convergence rate of $\mathcal{O}(T^{-1/2})$ in expectation. Moreover, if quantities $C_i(x_0),\,i\in \{2,3,4\}$ are known, then the coefficient associated with the rate-determining term can be lowered by a suitable choice of $\alpha$.
\end{enumerate}
}
\begin{proof}
    For statement $(i)$, since $\dfrac{1}{\sqrt{u+1}}$ is a positive and decreasing function for $u\ge0$, the following sum-integral inequalities hold for all $T\in\mathbb{N}$:
    \begin{gather}
        \sum_{t=0}^T \dfrac{1}{\sqrt{t+1}}\ge \int_0^{T+1}\dfrac{1}{\sqrt{u+1}}du=2(\sqrt{T+2}-1), \\
        \sum_{t=0}^T \dfrac{1}{t+1}\le 1+\int_{0}^{T}\dfrac{1}{u+1}du=1+\log(T+1).
    \end{gather}
    
    Then, the right-hand side of (\ref{convergence}) is upper bounded by
\begin{align}
    \label{order bound}
    \dfrac{1}{2\alpha(\sqrt{T+2}-1)}\left[M + 2\alpha C_1(x_0)R + 
    \frac{\alpha^2}{\mu}\left(\frac{C_2(x_0)}{2}+C_3(x_0)+R_*C_4(x_0)\right) (1+\log(T+1))  \right] = \mathcal{O}\left(\dfrac{\text{log}\,T}{\sqrt{T}}\right).
\end{align}

For statement $(ii)$, plugging in the proposed step size results in the right-hand side of (\ref{convergence}) becoming
\begin{equation}
    \dfrac{M}{\alpha\sqrt{T+1}}+ \dfrac{2C_1(x_0)R}{T+1}+\dfrac{\alpha\left(\frac{C_2(x_0)}{2}+C_3(x_0)+R_*C_4(x_0)\right)}{\mu\sqrt{T+1}}=\mathcal{O}(T^{-1/2})
\end{equation}

Note that, the above expression is a convex and continuously differentiable function of $\alpha$ on $(0,\infty)$, with its minimum value attained at the value where its derivative with respect to $\alpha$ is zero, assuming $\frac{C_2(x_0)}{2}+C_3(x_0)+R_*C_4(x_0) > 0$. Hence, if the quantities $C_i(x_0),\,i\in \{2,3,4\}$ are known, the bound can be improved with a value of $\alpha$ selected as 
\begin{equation}
    \alpha^\star=\sqrt{\dfrac{M\mu}{\frac{C_2(x_0)}{2}+C_3(x_0)+R_*C_4(x_0)}}.
\end{equation}
\end{proof}

\section{Solution to the Poisson Equation for Finite-State and Irreducible Markov Chains}

In this section, together with the preexisting conditions regarding the Markov process as presented in Assumption 3.1, we specialize our discussion of the Poisson equation to the domain of finite-state, irreducible Markov chains. It is well-known that a solution to the Poisson equation exists for such a chain, without having to explicitly assume that the chain is aperiodic (which together with our other assumptions would yield ergodicity). In light of the analysis in \cite{makowski1994poisson, kemeny1981generalization, bertsekas1976dynamic}, we first verify that a solution to the Poisson equation exists for our problem, defined in terms of the subgradient $g$ and the fundamental matrix of the chain, and, together with Assumption 3.3, show that the desired properties in Assumptions 3.4-3.5 are satisfied.  We also formulate a derivation for finite-state, irreducible, and periodic Markov chains, thereby explicitly showing that probabilistic mixing and convergence to the stationary distribution $\pi$ are not necessary for a solution to the Poisson equation to exist. Finally, we show the bounds obtained in terms of the fundamental matrix of the chain can be further upper bounded by quantities involving the mixing rate for aperiodic chains and the mixing rates of the $d$-step transition process for periodic chains where $d$ denotes the period, thereby recovering results in \cite{duchi2012ergodic} for $d=1$ and extending them to periodic chains.

Throughout this section, let $P\in\mathbb{R}^{|
\mathcal{X}|\times|\mathcal{X}|}$ denote the one-step transition matrix of the Markov chain and $\pi \in \mathbb{R}^{1\times|\mathcal{X}|}$ denote the stationary probability vector, which uniquely satisfies the set of equations $\pi=\pi P$ and $\pi \mathds{1}=1$, where $\mathds{1}$ denotes the column vector whose elements are all equal to 1. 

\subsection{General Irreducible Chains}

Here we present several results that are applicable to all finite-state and irreducible Markov chains. The analysis herein is largely a verification of existing results on the solution to the Poisson equation (see \cite{asmussen2003applied,makowski1994poisson, bertsekas1976dynamic}), adapted to our mirror descent problem.

\textbf{Lemma 5.1.1.} \emph{Let $W, U\in \mathbb{R}^{|\mathcal{X}|\times d}$ for some $d\in \mathbb{Z}^+$ such that $\pi U=0$. Then, the solution to the equation $(I-P)W=U$ subject to $\pi W = 0$ is $W=ZU$, where $Z\coloneqq(I-P+\mathds{1}\pi)^{-1}$ is the fundamental matrix of the Markov chain}. 
\begin{proof}
    See Appendix B.
\end{proof}

\textbf{Corollary 5.1.2.} \emph{Suppose $U=(I-\mathds{1}\pi)H$ for some $H\in \mathbb{R}^{|\mathcal{X}|\times d}$. Then, the solution to the equation $(I-P)W=U$ with $\pi W=0$ is $W=YH$, where $Y\coloneqq Z-\mathds{1}\pi$}.
\begin{proof}
    See Appendix B.
\end{proof}

\textbf{Theorem 5.1.3.} \emph{Assumption 3.4 is satisfied with the functions $B_1(x), B_2(x)$ defined as $B_i(x)\coloneqq  \\\sum_{x'\in \mathcal{X}} |Y_{x\,x'}| A_i(x'), \, \, i \in \{1,2 \}.$ This leads to Assumption 3.5 being automatically satisfied, with the pertinent quantities being expressed in terms of the maximum values of $A_1(x)$ and $A_2(x)$ across the finite state space, along with $\lVert Y\rVert_\infty$.}
\begin{proof}
    Here we include a sketch of the proof, with the full proof available in Appendix B.

    For our problem, the Poisson equation can be written as follows:
\begin{equation}
    V_\theta(x)-\sum_{x'\in \mathcal{X}}V_\theta(x')P_{x\,x'}= g(\theta,x)-\sum_{x'\in\mathcal{X}}\pi(x')g(\theta,x'), \, \, \forall \,x\in\mathcal{X}.
\end{equation}

Omitting a few steps, we have that for each state $x\in \mathcal{X}$, the solution to the Poisson equation is $V_\theta(x)=\sum_{x'\in \mathcal{X}} Y_{x\,x'}\,g(\theta,x')$ by Corollary 5.1.2. Hence, utilizing Assumption 3.3, we arrive at the following for all $x \in \mathcal{X}$ and $\theta, \theta' \in \widetilde{\Theta}$:
    \begin{gather}
        \lVert V_\theta(x)\rVert_*\le \sum_{x'\in \mathcal{X}}|Y_{x\,x'}| \, \lVert g(\theta,x')\rVert_* \le \sum_{x'\in \mathcal{X}}|Y_{x\,x'}| A_1(x')= B_1(x),\nonumber\\
        \lVert V_{\theta'}(x)-V_\theta(x)\rVert\le \lVert \theta'-\theta\rVert \sum_{x'\in \mathcal{X}}|Y_{x\,x'}|A_2(x') = B_2(x)\lVert \theta'-\theta\rVert.
    \end{gather}

    Hence, through $\lVert Y\rVert_\infty$ and the maximum values of $A_1(x)$ and $A_2(x)$ across $\mathcal{X}$, it follows that Assumption 3.5 is satisfied.
\end{proof}

The finite-time convergence rate derived in the previous section is therefore a function of $\lVert Y \rVert_\infty$, which is a property of the Markov chain itself. Note that, so far, we have not made any assumptions on whether the Markov chain at hand is periodic or aperiodic, since, for finite-state Markov chains, irreducibility alone is sufficient for the quantity $\pi$ to be uniquely defined and for the Poisson equation to have a unique solution $V$, where $\pi V=0$.
Now we can consider the cases of whether the Markov chain is periodic or aperiodic to derive an upper bound for $\lVert Y\rVert_\infty$.

\subsection{A Bound on \texorpdfstring{$\lVert Y \rVert_\infty$}{||Y||∞} for Periodic Chains}

\textbf{Definition:} A Markov chain is called \emph{periodic} with a period of $d \in \mathbb{Z}^+$ if the greatest common divisor of the elements of the set $\{n \in \mathbb{Z}^+: P^{(n)}_{ii}>0, i\in \{1, 2, ..., |\mathcal{X}|\} \}$ is equal to $d$. If $d=1$, the chain is denoted as being \emph{aperiodic}.

\textbf{Lemma 5.2.1 }(Theorem 1.8.4 of \cite{norris1998markov}). \emph{A Markov chain with a period of $d$ admits a partition of its states $C_0 \cup C_1 \cup \dots \cup C_{d-1} = \mathcal{X}$ such that the following hold (where the set indices are interpreted modulo $d$):
\begin{enumerate}[label=(\roman*)]
\item $P^{(n)}_{ij} >0$ only if $i \in C_r$ and $j \in C_{r+n}$ for some $r$,
\item $P^{(nd)}_{ij} >0$ for all sufficiently large $n$, for all $i,j \in C_r$, for all $r$.
\end{enumerate}}

\begin{proof}
    See Theorem 1.8.4. of \cite{norris1998markov}.
\end{proof}

This directly implies that, following a possible reordering of states, the transition matrix P can be written in the following cyclic block form:

\begin{equation}
    \label{periodic P}
    P=
    \begin{pmatrix}
    0      & P_0    & 0      & \cdots & 0 \\
    0      & 0      & P_1    & \cdots & 0 \\
    \vdots & \vdots & \vdots & \ddots & \vdots \\
    0      & 0      & 0      & \cdots & P_{d-2}\\
    P_{d-1}& 0      & 0      & \cdots & 0
    \end{pmatrix}
\end{equation}

where $P_r\in \mathbb{R}^{|C_r| \times |C_{r+1}|}$ for all $r$, assuming $C_{md+r}=C_r$ for all $m\in \mathbb{Z} $. Also, each row of $P_r$ sums to $1$, for all $r$.

Now, consider the Poisson equation of the form $(I-P)w=u$ where $w,u\in\mathbb{R}^{|\mathcal{X}|}$ and $\pi u=0$, subject to $\pi w=0$. Analyzing the properties of this equation with the cyclic block form of $P$ from (\ref{periodic P}) in mind will yield several useful results. 

Suppose the vectors $u$ and $w$ are partitioned as $u = [u_0^\top \, \, u_1^\top \, \, ... \, \, u_{d-1}^\top]^\top$, $w = [w_0^\top \, \, w_1^\top \, \, ... \, \, w_{d-1}^\top]^\top$ where $u_r, w_r\in \mathbb{R}^{|C_r|}$ for all $r$. Then, through the block decomposition in (\ref{periodic P}), the Poisson equation can be written in terms of each $w_r$ and $u_r$ as follows, where the indices are interpreted modulo $d$:
\begin{equation}
    \label{one step Poisson}
    w_r-P_r\,w_{r+1} = u_r.
\end{equation}
Let $Q_r\coloneqq \prod_{l=0}^{d-1}P_{r+l}$ and $S_r\coloneqq\sum_{l=0}^{d-1}\left[\prod_{k=1}^{l}P_{r+k-1}\right]u_{r+l}$, where indices are interpreted modulo $d$ and the product operator returns the identity matrix at $l=0$. Then, recursively applying (\ref{one step Poisson}) $d$ times yields
\begin{equation}
    \label{d step Poisson}
    (I-Q_r)w_r=S_r.
\end{equation}
Additionally, note that 
$P^d = \mathrm{diag}(Q_0, Q_1,...,Q_{d-1})$ is a block-diagonal matrix. Hence, through statement (ii) of Lemma 5.2.1, $Q_r$ is irreducible and aperiodic for all $r$. Since each $Q_r$ is irreducible, a unique stationary distribution $\mu_r$ exists, satisfying $\mu_r Q_r=\mu_r$ and $\mu_r \mathds{1} = 1$.

\textbf{Lemma 5.2.2.} \emph{Let the stationary row vector $\pi$ be partitioned as $\pi = [\pi_0 \, \, \pi_1 \, \, ... \, \, \pi_{d-1}]$ where $\pi_r\in \mathbb{R}^{|C_r|}$ for all $r$. Then, the stationary distribution $\mu_r$ of $Q_r$ is equal to $\mu_r=d\pi_r$ for all $r$.} 
\begin{proof}
    See Appendix B.
\end{proof}

Moreover, note that $\mu_rS_r= \sum_{l=0}^{d-1}\left[ \mu_r \prod_{k=1}^{l}P_{r+k-1}\right]u_{r+l}=\sum_{l=0}^{d-1} \mu_{r+l}u_{r+l}=d\pi u=0$, given that $\pi u = 0$ (indices are interpreted modulo $d$). Hence, (\ref{d step Poisson}) represents a valid Poisson equation.

We are now able to state the following theorem on a bound for $\lVert Y \rVert_\infty$ for periodic Markov chains:

\textbf{Theorem 5.2.3.} \emph{For an irreducible, finite-state Markov chain with period $d$, $\lVert Y \rVert_\infty \le 2((d-1)+d\max_r \lVert Y_r\rVert_\infty)$, where $Y_r\coloneqq (I-Q_r+\mathds{1}_{|C_r|}\mu_r)^{-1}-\mathds{1}_{|C_r|}\mu_r$ for each $r \in \{0, 1, ..., d-1\}$.}

\begin{proof}
    We have the following:
\begin{align}
    \label{Yu}
    \lVert Y\rVert_\infty &= \sup_{h\in \mathbb{R}^{|\mathcal{X}|}, \,\lVert h\rVert_\infty \le 1}  \lVert Yh\rVert_\infty \nonumber \\
    &= \sup_{h\in \mathbb{R}^{|\mathcal{X}|}, \,\lVert h\rVert_\infty \le 1}  \lVert Y(h-\mathds{1}\pi h)\rVert_\infty && \text{($Y\mathds{1}=Z\mathds{1}-\mathds{1} = 0$)} 
    \nonumber \\
    &\le \sup_{u\in \mathbb{R}^{|\mathcal{X}|}, \, \pi u = 0, \,\lVert u\rVert_\infty \le 1}  2\lVert Yu\rVert_\infty && \text{($\pi (h-\mathds{1}\pi h)=0$; $\lVert h - \mathds{1}\pi h\rVert_\infty \le 2 \lVert h \rVert_\infty\le 2$)} 
    \nonumber \\
\end{align}

   For every $u$ satisfying $\pi u=0$, $w=Yu$ is the unique solution that satisfies both the Poisson equation, $(I-P)w=u$, and the zero stationary average condition, $\pi w=0$. Suppose $w$ is partitioned as $w = [w_0^\top \, \, w_1^\top \, \, ... \, \, w_{d-1}^\top]^\top$ where $w_r\in \mathbb{R}^{|C_r|}$ for all $r$. Then, noting that $(I-Q_r)\mathds{1}_{|C_r|}=0$, we can modify (\ref{d step Poisson}), yielding $(I-Q_r)(w_r-\mathds{1}_{|C_r|}\mu_rw_r)=S_r$, where $\mu_r(w_r-\mathds{1}_{|C_r|}\mu_rw_r)=0$. Since the unique solution for $\bar{w}_r$ in the series of equations $\mu_r\bar{w}_r=0$ and $(I-Q_r)\bar{w}_r=S_r$ is $\bar{w}_r=Y_rS_r$, we have

    \begin{equation}
        w_r= \mathds{1}_{|C_r|}\mu_r w_r + Y_rS_r.
    \end{equation}

    Hence,

    \begin{equation}
        \label{w bound}
        \lVert w_r \rVert_\infty \le |\mu_rw_r|+\lVert Y_r \rVert_\infty \lVert S_r\rVert_\infty
    \end{equation}

Note that 

\begin{align}
    \label{S bound}
    \lVert S_r\rVert_\infty = \left \lVert \sum_{l=0}^{d-1}\left[\prod_{k=1}^{l}P_{r+k-1}\right]u_{r+l} \right \rVert_\infty
    \le \sum_{l=0}^{d-1} \left \lVert \prod_{k=1}^{l}P_{r+k-1}\right \rVert_\infty \lVert u_{r+l} \rVert_\infty = \sum_{l=0}^{d-1} \lVert u_{r+l} \rVert_\infty \le d \lVert u \rVert_\infty
\end{align}

To bound $|\mu_r w_r|$, we do the following (note that all indices are interpreted modulo $d$):

\begin{align}
    \label{uniform bound}
    (\ref{one step Poisson}), \text{Lemma 5.2.2} \implies \mu_rw_r- \mu_{r+1}w_{r+1}&=\mu_ru_r \,\text{ for all $r$} \nonumber\\
    \implies \mu_rw_r-\mu_{r+k}w_{r+k}&=\sum_{l=0}^{k-1}\mu_{r+l}u_{r+l} \,\text{ for all $r, k$} \nonumber \\
    \implies |\mu_rw_r - \mu_{r+k}w_{r+k} | &\le \sum_{l=0}^{k-1} \lVert \mu_{r+l}\rVert_1 \lVert u_{r+l}\rVert_\infty \nonumber \\
    &\le (d-1)\lVert u \rVert_\infty \,\text{for all $r$, for all $k \in \{0, 1, ..., d-1 \}$} \nonumber \\ \nonumber \\ \implies |\mu_rw_r-\mu_lw_l|&\le (d-1)\lVert u\rVert_\infty \, \text{for all $r, l$}.\nonumber \\
\end{align}

Note that $\sum_{r=0}^{d-1}\mu_rw_r=d\pi w=0$. Hence, $\min_{k\in \{0,1,...,d-1 \}} \mu_kw_k \le 0, \,\,\max_{k\in \{0,1,...,d-1 \}} \mu_kw_k \ge 0$. Then, setting $m\in \arg \max_{k\in \{0,1,...,d-1 \}} \mu_kw_k$, $n\in \arg \min_{k\in \{0,1,...,d-1 \}} \mu_kw_k$ in (\ref{uniform bound}) yields

\begin{align}
    \label{mu w bound}
    |\mu_rw_r|\le \max \{|\mu_mw_m|,|\mu_nw_n|\}&=\max \{\mu_mw_m,-\mu_nw_n\}\nonumber \\ &\le \mu_mw_m-\mu_nw_n=|\mu_mw_m-\mu_nw_n| \le (d-1)\lVert u\rVert_\infty \,\text{for all $r$}.
\end{align}

Combining (\ref{w bound}), (\ref{S bound}), and (\ref{mu w bound}), we get

\begin{align}
    \label{final wr bound}
    \lVert w_r \rVert_\infty \le [(d-1)+d\lVert Y_r\rVert_\infty] \lVert u \rVert_\infty \le [(d-1)+d \max_{k\in \{0, 1, ..., d-1 \}}\lVert Y_k\rVert_\infty] \lVert u \rVert_\infty \nonumber \\
    \implies \lVert w \rVert_\infty = \max_{r\in \{0, 1, ..., d-1 \}}\lVert w_r \rVert_\infty  \le [(d-1)+d \max_k\lVert Y_k\rVert_\infty] \lVert u \rVert_\infty
\end{align}

Finally, we conclude the proof as follows:

For any $u$ satisfying $\lVert u\rVert_\infty\le1$ and $\pi u=0$ and for the resulting $w$ that satisfies $w=Yu$, $\lVert Yu\rVert_\infty=\lVert w\rVert_\infty \le [(d-1)+d \max_{k\in \{0, 1, ..., d-1 \}}\lVert Y_k\rVert_\infty] \lVert u \rVert_\infty \le [(d-1)+d \max_{k\in \{0, 1, ..., d-1 \}}\lVert Y_k\rVert_\infty].$

Taking the supremum of the left-hand side over all vectors $u$ subject to $\lVert u\rVert_\infty\le1$ and $\pi u=0$,

\begin{equation}
    \lVert Y\rVert_\infty \le \sup_{u\in \mathbb{R}^{|\mathcal{X}|}, \, \pi u = 0, \,\lVert u\rVert_\infty \le 1}  2\lVert Yu\rVert_\infty \le 2((d-1)+d \max_{k\in \{0, 1, ..., d-1 \}}\lVert Y_k\rVert_\infty)
\end{equation}

\end{proof}

\subsection{The Special Case of Ergodic Markov Chains}

In the context of irreducible, finite-state, and aperiodic Markov chains, it is quite well-known that the following \emph{geometric ergodicity} property holds: There exist constants $C > 0$ and $\rho \in (0,1)$ such that, for all $t \in \mathbb{N}$ and for all states $x_0 \in \mathcal{X}$,
\begin{equation}
\label{geometric ergodicity}
    \lVert \,\mathbb{P}(X_t\in \cdot \,| X_0=x_0) - \pi \rVert_{TV} \le C\rho^t,
\end{equation}
where, for any two probability vectors $p$ and $q$ of the same dimension, $\lVert p-q\rVert_{TV} = \sum_{x\in \mathcal{X}}|p(x)-q(x) |$. (Note that we lift this definition from \cite{duchi2012ergodic}. Other definitions may insert a factor of $\frac{1}{2}$ in front of the expression on the right-hand side). 

In the context of the stochastic mirror descent algorithm over general ergodic processes, a suitable bound is found in \cite{duchi2012ergodic}: For our Markov chain, let $\tau_{\text{mix}}\coloneqq \max_{x \in \mathcal{X}} \inf\{t \in \mathbb{N},  \,\lVert \mathbb{P}(X_t\in \cdot |X_0=x) - \pi \rVert_{TV}\le T^{-1/2}\}$ where $T \in \mathbb{N}$ denotes the horizon. Note that this quantity satisfies the uniform mixing assumption (Assumption C) of \cite{duchi2012ergodic}, evaluated at $\epsilon=T^{-1/2}$. As a result, the suboptimality gap has an upper bound of $\mathcal{O}\left(\sqrt{\tau_{\text{mix}}/T}+\tau_{\mathrm{mix}}/{T}\right)$ in expectation, ignoring lower-order terms and constants that do not depend on the stochastic process \cite[Equation 3.4]{duchi2012ergodic}. In the case of a process that satisfies (\ref{geometric ergodicity}), a bound of $\,\,
\mathcal{O}\Biggl(\sqrt{\dfrac{\left[\log(C)+\frac{1}{2}\log(T)\right]_+}{\log(\rho^{-1}) \,T}+\dfrac{1}{T}}+\dfrac{\left[\log(C)+\frac{1}{2}\log(T)\right]_+}{\log(\rho^{-1}) \,T}+\dfrac{1}{T}\Biggr)\,\,
$ can be recovered using this method. 

The known-horizon result that we have obtained in Corollary 4.2, combined with Theorem 5.1.3, yields a bound of $\mathcal{O}\left(\sqrt{\dfrac{1+\lVert Y\rVert_\infty}{T}}+\dfrac{\lVert Y \rVert_\infty}{T}\right)$, ignoring constants that do not depend on the Markov chain at hand. This result on its own already represents a notable development, improving the convergence rate from $\mathcal{O}\left(\sqrt{(\log T)/{T}}\right)$ to $\mathcal{O}(T^{-1/2})$ for a fixed geometrically ergodic Markov chain, in terms of dependence on the horizon. In addition, as will be made clear in the next few lines, the quantity $\lVert Y\rVert_\infty$ can be used to obtain a bound depending on the mixing rate $\rho$ if the chain is aperiodic, or on the mixing rates $\rho_r, \,r\in \{0, 1, ..., d-1\}$ of the $d-$step transition chain if the chain is periodic with period $d$. This is significant as, owing to the presence of at least two of the eigenvalues of their transition matrices on the unit circle, periodic Markov chains do not satisfy a uniform ergodicity property in the way that it is needed and as such fall outside of the scope of stochastic processes considered in \cite{duchi2012ergodic}. The relation between $\lVert Y \rVert_\infty$ and mixing rates for periodic Markov chains will come through Theorem 5.2.3, where the diagonal blocks of the $d$-step transition matrix $P^d$ exhibit the ergodicity properties necessary to decompose $\lVert Y\rVert_\infty$ into terms $\lVert Y_k\rVert_\infty$ for $k\in \{0,1,...,d-1\}$, with each $k$th component behaving akin to a distinct aperiodic chain.

As previously defined, $Y=Z-\mathds{1}\pi=(I-P+\mathds{1}\pi)^{-1}-\mathds{1}\pi$. First, assume that the transition matrix $P$ belongs to an irreducible and aperiodic chain. Leveraging Theorems 2.4 and 3.1 of \cite{meyer1975role}, we remark that $Y$ equals $(I-P)^{\#}$, where the symbol $\#$ corresponds to a \emph{group inverse}, and subsequently satisfies $Y=\sum_{k=0}^{\infty}(P^k-\mathds{1}\pi)$. Then, utilizing (\ref{geometric ergodicity}), we get the following:
\begin{align}
    \lVert Y \rVert_\infty \le \sum_{k=0}^{\infty} \lVert P^k - \mathds{1}\pi \rVert_\infty = \sum_{k=0}^{\infty} \max_{x\in \mathcal{X}} \sum_{x'\in\mathcal{X}}|(P^k)_{x\,x'}-\pi_{x'}|
    \le \sum_{k=0}^{\infty} C \rho^k
    =\dfrac{C}{1-\rho}.
    \label{Y inf bound aperiodic}
\end{align}

Leveraging Theorem 5.2.3, we get the following bound for periodic chains:
\begin{align}
\label{Y inf bound periodic}
        \lVert Y\rVert_\infty \le 2((d-1)+d \max_{r\in \{0, 1, ..., d-1 \}}\lVert Y_r\rVert_\infty)\le 2\left(d-1+\frac{d\,C_{\text{max}}}{1-\rho_{\text{max}}}\right),
\end{align}

where $C^{(r)}>0$ and $\rho^{(r)} \in (0,1)$ satisfy $\sum_{x'\in C_r} |(Q_r^t)_{x\,x'}-(\mu_{r})_{x'}|\le C^{(r)}(\rho^{(r)})^t $ for all $x\in C_r$, all $t \in \mathbb{N}$ and all $r \in \{0, 1, ..., d-1\}$, and $C_{\max}\coloneqq \max_{r} C^{(r)}, \,\rho_{\max}\coloneqq \max_r \rho^{(r)}$.

Hence, leveraging (\ref{Y inf bound aperiodic}) and (\ref{Y inf bound periodic}), the suboptimality gap can be upper bounded by expressions involving the mixing rate for aperiodic chains, which is in line with what is implied by the analysis in \cite{duchi2012ergodic}, and can be extended to periodic chains using an expression involving the mixing rates of the $d-$step transition process. However, we would like to remark that the bound we have derived in terms of $\lVert Y\rVert_\infty$ can, in certain cases, represent a great improvement over mixing-time-type bounds, especially as the mixing parameter $\rho$ approaches $1$. We will illustrate this point through the following example, where in certain cases the bounds we have examined that depend on $(\log(\rho^{-1}))^{-1}$ or $(1-\rho)^{-1}$ can grow arbitrarily large as $\rho\rightarrow1^-$ whereas $\lVert Y\rVert_{\infty}$ remains bounded:

\textbf{Example 5.3.1.} Consider the two-state irreducible chain with the following transition matrix P, parameterized by $p\in (0,1]:$
\begin{equation}
    P =
\begin{pmatrix}
1-p & p \\
p & 1-p
\end{pmatrix}.
\end{equation}
By inspection, $\pi = [1/2 \,\,\,\, 1/2].$ We calculate $Y$ to be 
\begin{equation}
    Y = \frac{1}{4p}
\begin{pmatrix}
1 & -1 \\
-1 & 1
\end{pmatrix},
\end{equation}
yielding $\lVert Y\rVert_{\infty} = \dfrac{1}{2p}.$ Diagonalizing $P$, we compute $P^k-\mathds{1}\pi$ for $k\in\mathbb{N}$ as
\begin{equation}
    P^k-\mathds{1}\pi =\dfrac{(1-2p)^k}{2}
\begin{pmatrix}
1 & -1 \\
-1 & 1
\end{pmatrix}.
\end{equation}
Hence, $\max_{x\in\{1,2\}}\lVert \,\mathbb{P}(X_k\in \cdot \,| X_0=x) - \pi \rVert_{TV} = \lVert P^k-\mathds{1}\pi\rVert_\infty=|1-2p|^k$, with suitable choices for $C$ and $\rho$ being $C=1, \,\,\rho=|1-2p|$ (assuming $p\neq 1/2$).

Now consider a family of Markov chains with transition matrices $P(p)$, allowing us to examine the performance of our algorithm as $p$ varies. Notice that, as $p\rightarrow1^-$ (leading to $\rho \rightarrow 1^-$) and at any fixed $T>1$, $\lVert Y\rVert_\infty\rightarrow \frac{1}{2}^+$ and $(\log(|1-2p|^{-1}))^{-1}\rightarrow +\infty$, leading the mixing-time-dependent upper bound obtained through \cite{duchi2012ergodic} to diverge to infinity while ours remains bounded. Additionally, our bound allows for $p=1$ to be inserted, which corresponds to the two-state chain becoming periodic, while that of \cite{duchi2012ergodic} does not allow for such a case. We note that, as $p \rightarrow 0^+$, both approaches yield results that diverge to infinity, since at $p=0$ the transition matrix is equal to the identity matrix and therefore becomes reducible.

\section{Markov Chains with Countably Infinite State Spaces}

Proving the existence of a solution to the Poisson equation and verifying  that our assumptions are satisfied is more involved in the case of countable state-space chains. In particular, we leverage an  assumption regarding Lyapunov drift properties of the Markov chain, motivated by similar assumptions in \cite{meyn2012markov, douc2018markov,srikant2025rates}. Additionally, while it is well-established that solutions to the Poisson equation can be characterized using Lyapunov drift arguments (see \cite{glynn1996liapounov,herve2025computable,care2026stochastic,glynn2026computable}), we believe the assumptions we have laid out in this section work best with the specifics of our mirror descent problem, particularly for the verification of Assumptions 3.4 and 3.5. As we later show in more detail, we make use of a Lyapunov function $W(\cdot)$ and its square-root version, $W^{1/2}(\cdot)$, utilizing the latter to derive hitting-time bounds and to control first-order terms while using the former to control second-order terms involving the gradient and the solution to the Poisson equation.

\textbf{Assumption 6.1}.
The Markov chain $\{X_t\}$ is irreducible. Moreover, there exist a nonempty finite 
subset $F\subset\mathcal{X}$, a Lyapunov function $W: \mathcal{X}\rightarrow[1,\infty)$, constants $\lambda\in [0,1)$ and $K \in [1,\infty)$
that satisfy
\begin{equation}
    \label{drift}
    \mathbb{E}_x[W(X_1)]\le \lambda W(x)+ K \,\mathds{I}\{x\in F\} \, \,\, \forall \, x\in \mathcal{X}.
\end{equation}
Note that, through (\ref{drift}), Jensen's inequality, and the inequality $(a+b)^{1/2} \le a^{1/2}+b^{1/2}$ for $a,b\ge 0$, we get the following drift inequality for $W^{1/2}(\cdot)$, for all $x\in \mathcal{X}$:
\begin{align}
\label{drift sqrt}
    \mathbb{E}_x[W^{1/2}(X_1)]\le (\mathbb{E}_x[W(X_1)])^{1/2}\le (\lambda W(x)+K\mathds{I}\{x\in F\})^{1/2}\le \lambda^{1/2}W^{1/2}(x)+K^{1/2}\mathds{I}\{ x \in F\}.
\end{align}

\emph{Remark.} This assumption subsumes our results for the finite state-space case, where the set $F$ can be taken to be $\mathcal{X}$. For countable state-space chains, this condition continues to allow the chain to be periodic and hence, non-ergodic. However, we remark that this assumption might not be satisfied for certain classes of slowly mixing ergodic chains. 

\textbf{Lemma 6.2.} \emph{For any stopping time $\tau$ taking values in $\mathbb{N}^+ \cup \{ \infty\}$ and any $N \in \mathbb{N}^+$, the following inequality holds for all $x\in\mathcal{X}$:}
\begin{align}
    \mathbb{E}_x[N \land \tau]\le \mathbb{E}_x\left[\sum_{t=0}^{(N\land \tau)-1}W^{1/2}(X_t)\right]\le \dfrac{W^{1/2}(x)-1 +K^{1/2}\,\mathbb{E}_x\left[ \sum_{t=0}^{(N\land\tau)-1}\mathds{I}\{ X_t \in F\}\right]}{1-\lambda^{1/2}}.
\end{align}
\begin{proof}
    See Appendix C.
\end{proof}

\textbf{Corollary 6.3.} \emph{Let $\tau_F^{(1)} \coloneqq \inf \{t>0: X_t\in F \}$ and $\tau_F^{(k)}\coloneqq \inf \{t > \tau_F^{(k-1)}: X_t\in F \} \, \, \forall \, k>1$. Then, for all states $x \in \mathcal{X}$ such that $X_0=x$, $\tau_F^{(k)} < \infty$ almost surely for all $k \in \mathbb{N}^+$}.
\begin{proof}
    See Appendix C.
\end{proof}

\textbf{Lemma 6.4.} \emph{Owing to the irreducibility of the Markov chain, for a certain state $z \in F$, there exists a constant $L \in [1,\infty)$ that satisfies 
\begin{equation}
    \mathbb{E}_x\left[\sum_{t=0}^{\tau_z-1} \mathds{I}\{X_t\in F \}\right] \le L \, \,\, \forall \, x\in \mathcal{X},
\end{equation}
where $\tau_z \coloneqq \inf \{t>0: X_t=z \}$.}
\begin{proof}
    See Appendix C.
\end{proof}

\textbf{Theorem 6.5}. \emph{As a consequence of (\ref{drift sqrt}), and Lemma 6.2, and Lemma 6.4, the following are satisfied for all $x \in \mathcal{X}$:}
\begin{gather}
    (i) \,\mathbb{E}_x[\tau_z]\le \mathbb{E}_x\left[\sum_{t=0}^{\tau_z-1} W^{1/2}(X_t)\right ]\le \dfrac{W^{1/2}(x)+K^{1/2}L-1}{1-\lambda^{1/2}}, \nonumber \\ (ii) \, \dfrac{\mathbb{E}_x\left[ \sum_{t=0}^{\tau_z-1}W^{1/2}(X_t)\right]}{\mathbb{E}_x[\tau_z]} \le  \dfrac{W^{1/2}(x)+K^{1/2}-1}{1-\lambda^{1/2}}.
\end{gather}
\begin{proof}
    See Appendix C.
\end{proof}

We find it convenient to make the following assumption, which builds on Assumption 3.3:

\textbf{Assumption 6.6.} There exist constants $D_1, D_2 \in (0,\infty)$ such that the following hold for all $x \in \mathcal{X}$:
\begin{gather}
    A_1(x)\le D_1W^{1/2}(x),\nonumber\\
    A_2(x)\le D_2W^{1/2}(x).
\end{gather}

\textbf{Theorem 6.7.} \emph{Owing to Assumption 3.3, Theorem 6.5, Assumption 6.6 and \cite[Theorem 6.4]{makowski1994poisson}, there exists a valid solution to the Poisson equation for our problem as follows, for all $x\in\mathcal{X}$ and all $\theta\in \widetilde{\Theta}$:}
\begin{align}
\label{sol to pois eq hitting time}
    V_\theta(x)= \mathbb{E}_x\left[ \sum_{t=0}^{\tau_z-1}\bigl(g(\theta,X_t)-\bar{g}(\theta)\bigr)
     \right]=\mathbb{E}_x\left[ \sum_{t=0}^{\tau_z-1}g(\theta,X_t)
     \right]-\dfrac{\mathbb{E}_x[\tau_z]}{\mathbb{E}_z[\tau_z]}\,\mathbb{E}_z\left[ \sum_{t=0}^{\tau_z-1}g(\theta,X_t)
     \right]
\end{align}

\emph{Moreover, statements $(i)$ and $(ii)$ of Assumption 3.4 hold.} (Note: Suitable selections of $B_1(x), B_2(x)$ are provided in Appendix C.)

\begin{proof}
Here we outline a condensed version of the proof; the full proof is available in Appendix C.

The existence of the claimed Poisson solution follows from \cite[Theorem 6.4]{makowski1994poisson}, whose recurrence and integrability conditions can be verified using Theorem 6.5 and Assumption 6.6. 

By the hitting-time representation and Assumption 6.6,
we have (omitting several steps)
\begin{gather}
        \lVert V_\theta(x) \rVert_* \le D_1 \left\{\mathbb{E}_x\left[\sum_{t=0}^{\tau_z-1}  W^{1/2}(X_t)\right]+ \dfrac{\mathbb{E}_x[\tau_z]}{\mathbb{E}_z[\tau_z]}\,\mathbb{E}_z\left[ \sum_{t=0}^{\tau_z-1}W^{1/2}(X_t) \right] \right\},\nonumber\\
        \lVert V_{\theta'}(x)-V_\theta(x) \rVert\le  D_2\, \lVert \theta' - \theta\rVert \left\{\mathbb{E}_x\left[\sum_{t=0}^{\tau_z-1}  W^{1/2}(X_t)\right]+ \dfrac{\mathbb{E}_x[\tau_z]}{\mathbb{E}_z[\tau_z]}\,\mathbb{E}_z\left[ \sum_{t=0}^{\tau_z-1}W^{1/2}(X_t) \right] \right\}.
\end{gather}
Applying Theorem 6.5, we get suitable selections of $B_1(x)$ and $B_2(x)$, thus satisfying Assumption 3.4.
\end{proof}

\textbf{Theorem 6.8.} \emph{Through Theorem 6.7 and Assumptions 6.1 and 6.6, Assumption 3.5 holds.} (Note: Suitable selections of $C_1(x), C_2(x), C_3(x), C_4(x)$ are provided in Appendix C.)
\begin{proof}
    As before, we include a short description of the proof, with its full version available in Appendix C.

    Through Theorem 6.7, we observe that the quantities $\,\mathbb{E}_{x_0}[B_1(X_t)], \,\,\mathbb{E}_{x_0}[A_1(X_t)^2], \\\mathbb{E}_{x_0}[A_1(X_t)B_1(X_{t+1})],\,\,\mathbb{E}_{x_0}[A_1(X_t)B_2(X_{t+1})]$ are bounded from above by functions in terms of $\\\mathbb{E}_{x_0}[W^{1/2}(X_t)],\,\, \mathbb{E}_{x_0}[W(X_t)],\,$ and $\,\mathbb{E}_{x_0}[W^{1/2}(X_t)W^{1/2}(X_{t+1})]$.

Let $\{X_t\}_{t\in\mathbb{N}}$ be a sample path of the Markov chain, with initial state $X_0=x_0$ for some $x_0\in\mathcal{X}$. Applying recursion on the drift inequalities in Assumption 6.1 yields the following uniform bounds across $t\in\mathbb{N}$:
\begin{gather}
    \mathbb{E}_{x_0}[W(X_t)] \le  W(x_0)+\frac{K}{1-\lambda},\,\,\mathbb{E}_{x_0}[W^{1/2}(X_t)] \le W^{1/2}(x_0)+\frac{K^{1/2}}{1-\lambda^{1/2}}
\end{gather}

Bounding $\mathbb{E}_{x_0}[W^{1/2}(X_t)W^{1/2}(X_{t+1})]$ is straightforward, using the Cauchy-Schwarz inequality for nonnegative random variables:
\begin{equation}
    \mathbb{E}_{x_0}[W^{1/2}(X_t)W^{1/2}(X_{t+1})]\le \left( \mathbb{E}_{x_0}[W(X_t)]\,\mathbb{E}_{x_0}[W(X_{t+1})]\right)^{1/2} \le W(x_0)+\dfrac{K}{1-\lambda}.
\end{equation}

Hence, combining these bounds with those from Assumption 6.6 and Theorem 6.7, we have that Assumption 3.5 is satisfied, with selections of $C_i(x), i\in\{1,2,3,4 \}$ provided in Appendix C.
\end{proof}

\section{Conclusion}
In this paper, we have analyzed the MMD algorithm, a version of stochastic mirror descent where the objective is defined as the expectation of a quantity under the stationary distribution of a Markov chain and observations of the subgradient are corrupted by that chain along a given sample path $\{X_t \}_{t\in \mathbb{N}}$. We have shown that this algorithm yields a $\mathcal{O}(T^{-1/2})$ convergence rate in expectation when the horizon $T$ is known beforehand, which matches the rate of similar algorithms in stochastic (sub)gradient descent literature \cite{nemirovski2009robust}. A notable contribution that our analysis provides is that it incorporates the Poisson equation to explicitly verify that the algorithm converges even when the chain is not assumed to be ergodic, provided that several other assumptions are satisfied. Through this approach, we have improved the bounds derived in \cite{duchi2012ergodic}, both in terms of dependence on the horizon and in terms of the constants associated with the random process at hand. For irreducible and finite-state chains, we have derived an explicit bound on the expected value of the suboptimality gap for periodic chains. In the case of countably infinite state spaces, we invoke Lyapunov drift assumptions to obtain our results. For quantities involving the norm and Lipschitz properties of the subgradient and the solution to the Poisson equation, we have been able to derive upper bounds in expectation, conditioned on the initial state of the Markov chain. In particular, we allow for the unboundedness of these quantities over the state space of the chain, unlike what is assumed in \cite{paul2026stochastic}. Markov chains over more general, possibly uncountably infinite state spaces are an obvious direction to take should one wish to extend this work. Another possible direction is to examine the case where the transition kernel of the Markov chain is parameterized by the mirror descent iterate $\theta$, which would require generalizing the results of \cite{paul2026stochastic, care2026stochastic} in a way that is compatible with the possible unboundedness of quantities related to the subgradient and the solution to the Poisson equation and with certain classes of periodic chains.

\section*{Acknowledgments}

Research supported by NSF Grants CNS 23-12714, CCF 22-07547 and AFOSR Grant FA9550-24-1-0002.

We acknowledge the use of generative AI tools for literature review.

\bibliography{references}
\bibliographystyle{plainnat}

\newpage
\appendix
\section{Additional Details Related to the Graph Inference Problem}

\subsection{Proof of \texorpdfstring{
$\{(\omega_0)_i\}=\arg\min_{\omega_i\in\Omega_r}
\mathbb{E}_{Y^{(i)}\sim\pi_{Y^{(i)}}}
[f^{(i)}(\omega_i,Y^{(i)})]$
}{Uniqueness of the Minimizer}}

For any $\omega_i=[\hat{A}(i,:),\hat{b}(i)]\in\Omega_r$, we have $0\le f^{(i)}(\omega_i, Y^{(i)})\le -\log \left(\hat{b}(i)/2\right)\le f_{\text{max}}$ where $f_{\text{max}}\coloneqq-\log \left(b_{\text{min}}/2\right)$. Hence, we get the following:
\begin{align}
    \Delta^{(i)}(\omega_i,T)&\coloneqq\left| \frac{1}{T}\sum_{t=0}^{T-1} f^{(i)}({\omega_i},Y^{(i)}_t)-\mathbb{E}_{Y^{(i)}\sim\pi_{Y^{(i)}}}[f^{(i)}({\omega_i},Y^{(i)})] \right|\nonumber\\&=\left| \frac{1}{T}\sum_{t=0}^{T-1}\sum_{y\in\mathcal{Y}} \mathds{I}\{Y^{(i)}_t=y\}f^{(i)}(\omega_i, y)-\sum_{y\in\mathcal{Y}}\pi_{Y^{(i)}}(y)f^{(i)}(\omega_i,y) \right|\nonumber\\
    &=\left|\sum_{y\in\mathcal{Y}}f^{(i)}(\omega_i,y)\left(\frac{1}{T}\sum_{t=0}^{T-1} \mathds{I}\{Y^{(i)}_t=y \}-\pi_{Y^{(i)}}(y)\right) \right|\nonumber\\
    &\le f_{\text{max}} \lVert \hat{\pi}_{Y^{(i)}}^{(T)}-\pi_{Y^{(i)}}\rVert_1 
\end{align}
where $\hat{\pi}^{(T)}_{Y^{(i)}}(y) \coloneqq \frac{1}{T}\sum_{t=0}^{T-1}\mathds{I}\{Y^{(i)}_t=y\}$. From the ergodic theorem for Markov chains relating to time averages (which does not require aperiodicity), it follows that \(\lVert\hat{\pi}_{Y^{(i)}}^{(T)}-\pi_{Y^{(i)}} \rVert_1 \xrightarrow[T\to\infty]{\text{a.s.}}0\).  Hence, we have $ \sup_{\omega_i\in\Omega_r} \Delta^{(i)}(\omega_i,T) \xrightarrow[T\to\infty]{\text{a.s.}}0$.
Through this result, we have the following for any minimizer $\omega^\star_i\in\\ \arg\min_{\omega_i\in\Omega_r} \mathbb{E}_{Y^{(i)}\sim\pi_{Y^{(i)}}}[f^{(i)}(\omega_i,Y^{(i)})]$:
\begin{align}
    \mathbb{E}_{Y^{(i)}\sim\pi_{Y^{(i)}}}[f^{(i)}(\omega_i^\star,Y^{(i)})]&\le \mathbb{E}_{Y^{(i)}\sim\pi_{Y^{(i)}}}[f^{(i)}(({\omega_T^{\text{ML}}})_i,Y^{(i)})]\nonumber\\ &\le \frac{1}{T}\sum_{t=0}^{T-1} f^{(i)}(({\omega_T^{\text{ML}}})_i,{Y}^{(i)}_t)+\sup_{\omega_i\in\Omega_r} \Delta^{(i)}(\omega_i,T) \nonumber\\
    &\le \frac{1}{T}\sum_{t=0}^{T-1} f^{(i)}(\omega_i^\star,Y^{(i)}_t)+\sup_{\omega_i\in\Omega_r} \Delta^{(i)}(\omega_i,T) \nonumber\\
    &\le  \mathbb{E}_{Y^{(i)}\sim\pi_{Y^{(i)}}}[f^{(i)}(\omega_i^\star,Y^{(i)})]+ 2\sup_{\omega_i\in\Omega_r} \Delta^{(i)}(\omega_i,T)
\end{align}

Taking limits as $T\rightarrow\infty$ on a sample path along which $\omega^{\mathrm{ML}}_T \xrightarrow[T\to\infty]{} \omega_0$ and $ \sup_{\omega_i\in\Omega_r} \Delta^{(i)}(\omega_i,T)\\ \xrightarrow[T\to\infty]{}0$ and utilizing continuity, we get
\begin{gather}
    \mathbb{E}_{Y^{(i)}\sim\pi_{Y^{(i)}}}[f^{(i)}((\omega_0)_i,Y^{(i)})]= \mathbb{E}_{Y^{(i)}\sim\pi_{Y^{(i)}}}[f^{(i)}(\omega^\star_i,Y^{(i)})]\nonumber\\\implies (\omega_0)_i\in \arg\min_{\omega_i\in\Omega_r} \mathbb{E}_{Y^{(i)}\sim\pi_{Y^{(i)}}}[f^{(i)}(\omega_i,Y^{(i)})].
\end{gather}
All that remains is to show that $\arg\min_{\omega_i\in\Omega_r} \mathbb{E}_{Y^{(i)}\sim\pi_{Y^{(i)}}}[f^{(i)}(\omega_i,Y^{(i)})] = \{(\omega_0)_i \}$, i.e. that $(\omega_0)_i$ is the unique minimizer. To do that, we utilize properties of the Kullback-Leibler divergence, as shown in the next few lines. For any  $\omega_i\in\Omega_r$, the following holds:
\begin{align}
     \mathbb{E}_{Y^{(i)}\sim\pi_{Y^{(i)}}}[f^{(i)}(\omega_i,Y^{(i)})] -& \mathbb{E}_{Y^{(i)}\sim\pi_{Y^{(i)}}}[f^{(i)}((\omega_0)_i,Y^{(i)})]\nonumber\\&= \sum_{y\in\mathcal{Y}}\pi_{Y^{(i)}}(y)\log\dfrac{ P^{(i)}_{(\omega_0)_i}(y(p+1)\,|\,y(1:p))}{P^{(i)}_{\omega_i}(y(p+1)\,|\,y(1:p))}\nonumber \\
     &= \sum_{x\in\mathcal{X}}\sum_{z\in\{0,1\}}\pi_X(x)P_{(\omega_0)_i}^{(i)}(z|x)\log\dfrac{ P_{(\omega_0)_i}^{(i)}(z|x)}{ P_{\omega_i}^{(i)}(z|x)}\nonumber\\
&=\sum_{x\in\mathcal{X}}\underbrace{\pi_X(x)}_{>\,0}\underbrace{\sum_{z\in\{0,1\}}P_{(\omega_0)_i}^{(i)}(z|x)\log \dfrac{ P_{(\omega_0)_i}^{(i)}(z|x)}{ P_{\omega_i}^{(i)}(z|x)}}_{\ge \,0}\ge 0
\end{align}
with equality if and only if $P_{(\omega_0)_i}^{(i)}(z|x)=P_{\omega_i}^{(i)}(z|x)$ for all $x\in\mathcal{X}$ and $z\in\{0,1\}$. We will complete the proof by showing that, for any $\omega_1^{(i)},\omega_2^{(i)}\in\Omega_r$ such that $\omega_1^{(i)}\neq \omega_2^{(i)}$, there exists an $x\in\mathcal{X}$ and $z\in\{0,1\}$ such that $P_{\omega_1^{(i)}}^{(i)}(z|x)\neq P_{\omega^{(i)}_2}^{(i)}(z|x)$. Suppose that $\omega_1^{(i)}(p+1)\neq \omega_2^{(i)}(p+1)$. Then, choosing $x=0_p$ and $z=1$ yields $P_{\omega_1^{(i)}}^{(i)}(z|x)=\omega_1^{(i)}(p+1)/2\neq \omega_2^{(i)}(p+1)/2=P_{\omega^{(i)}_2}^{(i)}(z|x)$. Suppose instead that $\omega_1^{(i)}(p+1)= \omega_2^{(i)}(p+1)$ but $\omega_1^{(i)}(j)\neq \omega_2^{(i)}(j)$ for some $1\le j\le p$. Then, choosing $z=1$ and $x=e_j$ where $e_j$ is the $j$th column of the $p$-dimensional identity matrix yields $P_{\omega_1^{(i)}}^{(i)}(z|x)=\omega_1^{(i)}(j)+\omega_1^{(i)}(p+1)/2\neq \omega_2^{(i)}(j)+\omega_1^{(i)}(p+1)/2=P_{\omega^{(i)}_2}^{(i)}(z|x)$. Hence, we have that $(\omega_0)_i$ is the unique minimizer of $\mathbb{E}_{Y^{(i)}\sim\pi_{Y^{(i)}}}[f^{(i)}(\omega_i,Y^{(i)})]$.

Note that this proof establishes that $(\omega_0)_i$ is the unique minimizer of $\mathbb{E}_{Y^{(i)}\sim\pi_{Y^{(i)}}}[f^{(i)}(\omega_i,Y^{(i)})]$, without having to show its relation to the ML estimator (which yielded that $(\omega_0)_i$ is a valid minimizer). However, we believe that derivation is still useful in that it bridges the gap between previous work on this particular example problem and the broader setting of mirror descent explored in this paper, and also provides insight and motivation into deriving a suitable objective to be optimized with the mirror descent algorithm as the pointwise limit of an objective used in an ML estimation problem.

\subsection{Verification of Assumptions and Convergence}
In line with the framework we introduce in Section 3, our goal is to minimize $\bar{f}^{(i)}_M(u)$ where $u \in\Delta$, using observations of the stochastic gradient $g^{(i)}_M(u,Y^{(i)}_t)$ (which we assume are made available to the user at every time step) over a sample path $\{Y^{(i)}_t\}$ .
We let $h(u)=\sum_{j=1}^{p+1}u_j\log u_j$ for all $u\in\Delta$ with the convention that $0\log0\equiv0$. The Bregman divergence induced by $h$ then becomes the Kullback-Leibler divergence: $D_h(u,v)=D_{\text{KL}}(u\|v)=\sum_{j=1}^{p+1}u_j\log(u_j/v_j)\ge \frac{1}{2}\lVert u-v \rVert_1^2$ for all $u\in\Delta$ and $v\in \operatorname{relint}(\Delta)$, with strong convexity parameter $\mu=1$. We equip this algorithm with the $\lVert \cdot\rVert_1$ norm as the primal norm and $\lVert \cdot\rVert_\infty$ as the dual norm, with associated radii $R=2$ and $R_*=1$, respectively. We let $\{u^{(i)}_t\}_{t\in\mathbb{N}}$ be the iterates of our algorithm, with $u^{(i)}_0 = (1/(p+1))\mathds{1}$, yielding $D_h(u,u^{(i)}_0)\le \log (p+1)$ for all $u\in\Delta$. We update $u^{(i)}_{t+1}$ from $u^{(i)}_t$ as follows, for a suitably chosen step size $\epsilon_t$:
\begin{gather}
    u^{(i)}_{t+1}= \arg\min_{u \in \Delta}\, \langle g_M^{(i)}(u^{(i)}_t, Y^{(i)}_t), u \rangle + \frac{1}{\epsilon_t}D_h(u, u^{(i)}_{t})\nonumber\\
    \implies u^{(i)}_{t+1}(j)=\dfrac{u^{(i)}_{t}(j)\exp(-\epsilon_t\,g_M^{(i)}(u^{(i)}_t,Y^{(i)}_t)(j))}{\sum_{k=1}^{p+1} u^{(i)}_{t}(k)\exp(-\epsilon_t\,g_M^{(i)}(u^{(i)}_t,Y^{(i)}_t)(k))}\,\,\,\forall\, 1\le j\le p+1.
\end{gather}

Note that this implies $u^{(i)}_t\in \operatorname{relint}(\Delta)\,\,\forall\,t\in\mathbb{N}$, meaning the set defined in Assumption 3.2 (iv) can be taken to be $\operatorname{relint}(\Delta)$. We show that the stochastic gradient is Lipschitz in the $\lVert \cdot\rVert_1$ norm and that its $\lVert \cdot\rVert_\infty$ norm is bounded as follows:
\begin{align}
    g_M^{(i)}(u,y)&=-y(p+1)\dfrac{(1-b_{\min})[y(1:p);1/2]}{(1-b_{\min})\langle u(1:p), \, y(1:p) \rangle+((1-b_{\min})u(p+1)+b_{\min})(1/2)}\nonumber\\
    &+(1-y(p+1))\dfrac{(1-b_{\min})[y(1:p);1/2]}{1-(1-b_{\min})\langle u(1:p), \, y(1:p) \rangle-((1-b_{\min})u(p+1)+b_{\min})(1/2)}
\end{align}
\begin{gather}
\implies \lVert g_M^{(i)}(u,y) \rVert_\infty \le A_1,\,\,\lVert g_M^{(i)}(u,y)-g_M^{(i)}(u',y)\rVert_1\le A_2\lVert u-u' \rVert_1 \,\,\,\forall\,u,u'\in\operatorname{relint}(\Delta),\, y\in\mathcal{Y} \nonumber\\
\text{where } A_1\coloneqq\dfrac{2(1-b_{\min})}{b_{\min}},\,\,A_2\coloneqq \dfrac{4(1-b_{\min})^2(p+1/2)}{b_{\min}^2}
\end{gather}

We have so far shown that Assumptions 3.1-3.3 are satisfied for this problem. Utilizing Theorem 5.1.3, we have that Assumptions 3.4 and 3.5 are also satisfied. Hence, setting $(u^{(i)})^\star\coloneqq \Phi((\omega_0)_i)$, $\bar{u}^{(i)}_T\coloneqq\dfrac{\sum_{t=0}^T \epsilon_t u^{(i)}_t}{\sum_{t=0}^T \epsilon_t}$, invoking Theorem 4.1 and its corollary, we have that $\mathbb{E}[\bar{f}_M^{(i)}(\bar{u}^{(i)}_T)]-\bar{f}_M^{(i)}((u^{(i)})^\star)=\tilde{\mathcal{O}}(T^{-1/2})$, which together with Markov's inequality yields that $\bar{f}_M^{(i)}(\bar{u}^{(i)}_T)\xrightarrow[T\to\infty]{p}\bar{f}_M^{(i)}((u^{(i)})^\star)$. 

To show that $\bar{u}^{(i)}_T\xrightarrow[T\to\infty]{p}(u^{(i)})^\star$, we make the following argument: For any $\epsilon>0$, define the set $\mathcal{S}_\epsilon\coloneqq \{u^{(i)}\in\Delta\,|\, \lVert u^{(i)}-(u^{(i)})^\star \rVert_1\ge \epsilon \}$. For every $\epsilon>0$ such that $\mathcal{S}_\epsilon$ is nonempty, $\delta_\epsilon\coloneqq \min_{u^{(i)}\in\mathcal{S}_\epsilon} \bar{f}^{(i)}_M(u^{(i)})-\bar{f}^{(i)}_M((u^{(i)})^\star)$. Since $\bar{f}^{(i)}_M$ is continuous and the set over which the minimization takes place is compact, we have that $\delta_\epsilon$ is well-defined by the extreme value theorem. Moreover, since the minimizer of $\bar{f}^{(i)}_M$ is unique and is equal to $(u^{(i)})^\star$, we have that $\delta_\epsilon>0$. Hence, for every $\epsilon > 0$ such that $\delta_\epsilon$ is defined, we have that $\mathbb{P}(\lVert \bar{u}^{(i)}_T- (u^{(i)})^\star\rVert_1\ge \epsilon)\le \mathbb{P}(\bar{f}^{(i)}_M(\bar{u}_T^{(i)})-\bar{f}^{(i)}_M((u^{(i)})^\star)\ge\delta_\epsilon)=\tilde{\mathcal{O}}(T^{-1/2})$. (For $\epsilon>0$ such that $\mathcal{S}_\epsilon$ is empty, $\mathbb{P}(\lVert \bar{u}^{(i)}_T-(u^{(i)})^\star\rVert_1\ge\epsilon)=0$.)  Therefore, $\bar{u}^{(i)}_T\xrightarrow[T\to\infty]{p}(u^{(i)})^\star$ for all $i \in \{1,2,...,p\}$.

\section{Proofs Omitted in Section 5}

\textbf{Lemma 5.1.1.} \emph{Let $W, U\in \mathbb{R}^{|\mathcal{X}|\times d}$ for some $d\in \mathbb{Z}^+$ such that $\pi U=0$. Then, the solution to the equation $(I-P)W=U$ subject to $\pi W = 0$ is $W=ZU$, where $Z\coloneqq(I-P+\mathds{1}\pi)^{-1}$ is the fundamental matrix of the Markov chain}. 
\begin{proof}
    Since the Markov chain is irreducible, the null space of the matrix $I-P$ is simply $\text{span}(\mathds{1})$. Hence, for any two solutions $W_1, W_2$ to the equation $(I-P)W=U$, there must exist a vector $y\in \mathbb{R}^{d}$ such that $W_2=W_1 + \mathds{1}y^\top.$ Therefore, to ensure uniqueness of the solution, we may impose $\pi W=0$. Combining $(I-P)W=U$ and $\pi W=0$, we get $(I-P+\mathds{1}\pi)W=U$. The following is a concise proof of the invertibility of $(I-P+\mathds{1}\pi)$ for irreducible chains, of which a more general version is available in \cite{kemeny1981generalization}:

    Let $(I-P+\mathds{1}\pi)c=0$, or equivalently $(I-P)c=-\mathds{1}\pi c$. Multiplying from the left by $\pi$ and leveraging the fact that $\pi = \pi P$, we get $0=-\pi c$. Inserting this into the original equation, we get $(I-P)c=0$, or equivalently $Pc=c$. Since our Markov chain is irreducible, the only set of right eigenvectors corresponding to an eigenvalue of 1 is span($\mathds{1}$). Combining the equations $0 = \pi c$ and $c\in\text{span}(\mathds{1})$, we get $c=0$. Hence, $(I-P+\mathds{1}\pi)$ is invertible.

    Note that $\pi(I-P+\mathds{1}\pi)=\pi$, yielding $ \pi(I-P+\mathds{1}\pi)^{-1}U=\pi U=0.$ Additionally, $(I-P)(I-P+\mathds{1}\pi)^{-1}U=[(I-P+\mathds{1}\pi)-\mathds{1}\pi](I-P+\mathds{1}\pi)^{-1}U=U-\mathds{1}\pi U = U$.

    Hence, the solution to the equation $(I-P)W=U$ subject to $\pi W=0$ is $W=(I-P+\mathds{1}\pi)^{-1}U=ZU$.
\end{proof}

\textbf{Corollary 5.1.2.} \emph{Suppose $U=(I-\mathds{1}\pi)H$ for some $H\in \mathbb{R}^{|\mathcal{X}|\times d}$. Then, the solution to the equation $(I-P)W=U$ with $\pi W=0$ is $W=YH$, where $Y\coloneqq Z-\mathds{1}\pi$}.

\begin{proof}
    Note that $(I-P+\mathds{1}\pi)\mathds{1}=\mathds{1}$, so $Z\mathds{1}=\mathds{1}$. Hence, $ZU=Z(I-\mathds{1}\pi)H=(Z-\mathds{1}\pi)H=YH$.
\end{proof}

\textbf{Theorem 5.1.3.} \emph{Assumption 3.4 is satisfied with the functions $B_1(x), B_2(x)$ defined as $B_i(x)\coloneqq  \\\sum_{x'\in \mathcal{X}} |Y_{x\,x'}| A_i(x'), \, \, i \in \{1,2 \}.$ This leads to Assumption 3.5 being automatically satisfied, with the pertinent quantities being expressed in terms of the maximum values of $A_1(x)$ and $A_2(x)$ across the finite state space, along with $\lVert Y\rVert_\infty$.}
\begin{proof}
For our problem, the Poisson equation can be written as follows:
\begin{equation}
    V_\theta(x)-\sum_{x'\in \mathcal{X}}V_\theta(x')P_{x\,x'}= g(\theta,x)-\sum_{x'\in\mathcal{X}}\pi(x')g(\theta,x'), \, \, \forall \,x\in\mathcal{X}.
\end{equation}

Assuming the states of the Markov chain are numbered from 1 to $|\mathcal{X}|$, we express the Poisson equation in matrix form as follows:
\begin{equation}
    \label{matrixpoissoneq}
    (I-P)\mathcal{V}_\theta= (I-\mathds{1}\pi)\mathcal{G}_\theta
\end{equation}
where $\mathcal{V}_\theta$ and $\mathcal{G}_\theta$ are both $|\mathcal{X}|\times n$ matrices satisfying $\mathcal{V_\theta}(i,\cdot)=V_\theta(i)^\top$, $\mathcal{G_\theta}(i,\cdot)=g(\theta,i)^\top$, for all $i \in \{1, 2, ..., |\mathcal{X}|\}$.

    By Corollary 5.1.2, the solution to the matrix Poisson equation in (\ref{matrixpoissoneq}) satisfying $\pi \mathcal{V}_\theta=0$ is $\mathcal{V}_\theta=Y\mathcal{G}_\theta$ (note that this matrix form of the solution of the Poisson equation also appears in \cite{asmussen2003applied,makowski1994poisson, bertsekas1976dynamic}). Then, for each state $x\in \mathcal{X}$ the solution to the Poisson equation is $V_\theta(x)=\sum_{x'\in \mathcal{X}} Y_{x\,x'}\,g(\theta,x')$. Hence, utilizing Assumption 3.3, we arrive at the following for all $x \in \mathcal{X}$ and $\theta, \theta' \in \widetilde{\Theta}$:
    
    \begin{equation*}
        \lVert V_\theta(x)\rVert_*\le \sum_{x'\in \mathcal{X}}|Y_{x\,x'}| \, \lVert g(\theta,x')\rVert_* \le \sum_{x'\in \mathcal{X}}|Y_{x\,x'}| A_1(x')= B_1(x),
    \end{equation*}  
    \begin{align}
        \lVert V_{\theta'}(x)-V_\theta(x)\rVert&\le\left\lVert \sum_{x'\in \mathcal{X}}Y_{x\,x'}(\tilde{g}(\theta',x')+r(\theta')-\tilde{g}(\theta,x')-r(\theta)) \right\rVert \nonumber\\
            &=\left\lVert \sum_{x'\in \mathcal{X}}Y_{x\,x'}(\tilde{g}(\theta',x')-\tilde{g}(\theta,x')) \right\rVert && \text{($Y\mathds{1}=0$)} \nonumber\\
            &\le \sum_{x'\in \mathcal{X}}|Y_{x\,x'}| \, \lVert \tilde{g}(\theta',x')-\tilde{g}(\theta,x')\rVert \nonumber\\
            &\le \lVert \theta'-\theta\rVert \sum_{x'\in \mathcal{X}}|Y_{x\,x'}|A_2(x') = B_2(x)\lVert \theta'-\theta\rVert.
    \end{align}

    Let $\Tilde{A}_i \coloneqq \max_{x\in\mathcal{X}}A_i(x), \, \, i\in \{1,2\} .$ Let $C_1 \coloneqq \Tilde{A}_1\lVert Y \rVert_\infty$, $C_2 \coloneqq \Tilde{A}_1^2$, $C_3 \coloneqq \Tilde{A}_1^2\lVert Y \rVert_\infty$, and $C_4 \coloneqq \Tilde{A}_1\Tilde{A}_2\lVert Y \rVert_\infty$. Then, it is clear that, for any sample path $\{X_t\}$, $B_1(X_t)\le C_1$, $A_1(X_t)^2\le C_2$,\\ $A_1(X_t)B_1(X_{t+1})\le C_3$, and $A_1(X_t)B_2(X_{t+1}) \le C_4$ for all $t \in \mathbb{N}$. This automatically satisfies Assumption 3.5.
    \end{proof}

\textbf{Lemma 5.2.2.} \emph{Let the stationary row vector $\pi$ be partitioned as $\pi = [\pi_0 \, \, \pi_1 \, \, ... \, \, \pi_{d-1}]$ where $\pi_r\in \mathbb{R}^{|C_r|}$ for all $r$. Then, the stationary distribution $\mu_r$ of $Q_r$ is equal to $\mu_r=d\pi_r$ for all $r$.} 

\begin{proof}
    Noting that $\pi = \pi P$ and using the block decomposition of $P$ from (\ref{periodic P}), we get the following for each $r$ (where $r$ is interpreted modulo $d$):

    \begin{equation}
        \label{periodic pi}
        \pi_r =\pi_{r-1}P_{r-1}
    \end{equation}

    Recursively applying (\ref{periodic pi}) $d$ times and interpreting the indices modulo $d$, we get

    \begin{equation}
        \pi_r=\pi_{r-d}P_{r-d}P_{r-d+1}...P_{r-1}=\pi_rQ_r
    \end{equation}

    Hence, $\mu_r$ must be equal to a potentially scaled version of $\pi_r$, for all $r$. 

    As a consequence of (\ref{periodic pi}), $\pi_r \mathds{1}_{|C_r|}=\pi_{r-1} P_{r-1}\mathds{1}_{|C_r|}=\pi_{r-1}\mathds{1}_{|C_{r-1}|} \implies \pi_r \mathds{1}_{|C_r|} = \pi_l \mathds{1}_{|C_{l}|}$ for all $r$ and $l$ (interpreted modulo $d$).

    Let $c \coloneqq \pi_0 \mathds{1}_{|C_{0}|} = \pi_1 \mathds{1}_{|C_{1}|}= \,... = \pi_{d-1} \mathds{1}_{|C_{d-1}|}$. Then,

    \[
    1=\pi \mathds{1}=\sum_{r=0}^{d-1}\pi_r\mathds{1}_{|C_r|}=dc \implies c = \dfrac{1}{d}
    \]

    As a result, $\mu_r=d\pi_r$ for all $r$.
\end{proof}

\section{Proofs Omitted in Section 6}

\textbf{Lemma 6.2.} \emph{For any stopping time $\tau$ taking values in $\mathbb{N}^+ \cup \{ \infty\}$ and any $N \in \mathbb{N}^+$, the following inequality holds:}
\begin{equation}
  \mathbb{E}_x[N \land \tau]\le \mathbb{E}_x\left[\sum_{t=0}^{(N\land \tau)-1}W^{1/2}(X_t)\right] \le \dfrac{W^{1/2}(x)-1 +K^{1/2}\,\mathbb{E}_x\left[ \sum_{t=0}^{(N\land\tau)-1}\mathds{I}\{ X_t \in F\}\right]}{1-\lambda^{1/2}} \, \, \,\forall \, x \in \mathcal{X}   
\end{equation}

\begin{proof}
\begin{align}
    W^{1/2}(X_{N\land\tau})-W^{1/2}(X_0)&=\sum_{t=0}^{(N\land \tau)-1}\bigl(W^{1/2}(X_{t+1})-W^{1/2}(X_t)\bigr)\nonumber\\ &= \sum_{t=0}^{N-1} \bigl(W^{1/2}(X_{t+1})-W^{1/2}(X_t)\bigr)\, \mathds{I}\{t < \tau\}. 
\end{align}

Note that $W^{1/2}(X_{N\land\tau})\ge 1$. Additionally, since $\tau$ is a stopping time, $\mathds{I}\{t < \tau \}$ is $\mathcal{F}_t$-measurable. Taking expectation under the distribution $\mathbb{P}_x$ and utilizing (\ref{drift sqrt}), we get the following: 
\begin{align}
    1&\le W^{1/2}(x)+ \mathbb{E}_{x}\left[ \sum_{t=0}^{N-1} (W^{1/2}(X_{t+1})-W^{1/2}(X_t))\, \mathds{I}\{t < \tau \} \right] \nonumber \\
    &=W^{1/2}(x)+\mathbb{E}_{x}\left[  \sum_{t=0}^{N-1}  \mathbb{E}_{x}\left[(W^{1/2}(X_{t+1})-W^{1/2}(X_t))\, \mathds{I}\{t < \tau \}  \Big| \mathcal{F}_t\right]  \right] \nonumber \\
    &= W^{1/2}(x)+ \mathbb{E}_{x}\left[ \sum_{t=0}^{N-1} \left(\mathbb{E}_{x}[W^{1/2}(X_{t+1})|X_t]-W^{1/2}(X_t)\right) \mathds{I}\{t < \tau\}     \right] \nonumber \\
    &\le W^{1/2}(x)+\mathbb{E}_{x}\left[\sum_{t=0}^{N-1}((\lambda^{1/2}-1)W^{1/2}(X_t)+K^{1/2}\,\mathds{I}\{X_t\in F\})\mathds{I}\{t<\tau\} \right] \nonumber \\
    &\leq W^{1/2}(x) - (1-\lambda^{1/2})\,\mathbb{E}_{x}\left[\sum_{t=0}^{(N\land \tau)-1}W^{1/2}(X_t)  \right] + K^{1/2}\,\mathbb{E}_{x}\left[ \sum_{t=0}^{(N\land \tau)-1}\mathds{I}\{X_t\in F\}  \right].
\end{align}

Rearranging the last inequality and noting that $W^{1/2}(x) \ge 1 \, \,\forall \, x \in \mathcal{X}$, we get the claimed statement.
\end{proof}

\textbf{Corollary 6.3.} \emph{Let $\tau_F^{(1)} \coloneqq \inf \{t>0: X_t\in F \}$ and $\tau_F^{(k)}\coloneqq \inf \{t > \tau_F^{(k-1)}: X_t\in F \} \, \, \forall \, k>1$. Then, for all states $x \in \mathcal{X}$ such that $X_0=x$, $\tau_F^{(k)} < \infty$ almost surely for all $k \in \mathbb{N}^+$}.
\begin{proof}
    For every sample path with $X_0=x$ and for every $k \in \mathbb{N}^+$, the following holds:
    \begin{equation}
\sum_{t=0}^{(N\land\tau_F^{(k)})-1}\mathds{I}\{ X_t \in F\} \le k-1 + \mathds{I}\{x \in F \}
    \end{equation}
Utilizing Lemma 6.2 and the Monotone Convergence Theorem, we get
\begin{equation}
    \mathbb{E}_x[\tau_F^{(k)}]\le \dfrac{W^{1/2}(x)-1 + K^{1/2}(k-1+\mathds{I}\{x \in F \})}{1-\lambda^{1/2}} < \infty \, \, \,\forall \, x \in \mathcal{X}, \, \forall\,k\in \mathbb{N}^+.
\end{equation}
Hence, $\mathbb{P}_x(\tau_F^{(k)}=\infty)=0$, i.e. given any initial state $X_0=x$, $\tau_F^{(k)}$ is almost surely finite, for all $k \in \mathbb{N}^+$.
\end{proof}

\textbf{Lemma 6.4.} \emph{Owing to the irreducibility of the Markov chain, for a certain state $z \in F$, there exists a constant $L \in [1,\infty)$ that satisfies 
\begin{equation}
    \mathbb{E}_x\left[\sum_{t=0}^{\tau_z-1} \mathds{I}\{X_t\in F \}\right] \le L \, \,\, \forall \, x\in \mathcal{X},
\end{equation}
where $\tau_z \coloneqq \inf \{t>0: X_t=z \}$.}
\begin{proof}
    Noting that $\tau_F^{(k)}$ is almost surely finite for all $k\in \mathbb{N}^+$ for any initial state $X_0 = x$, consider the following random process $\{Y_n \}_{n\in\mathbb{N}}$: $Y_n \coloneqq X_{\tau_F^{(n+1)}} \, \, \forall \, n \in \mathbb{N}$, similar to the argument in the proof of \cite[Lemma 7.1.2]{bremaud2020markov}. This process is a Markov chain over the finite set $F$, with transition probabilities $\mathbb{P}(Y_1=j|Y_0=i)=\mathbb{P}_i(X_{\tau_F^{(1)}}=j) \, \,\forall \, i,j\in F$. Moreover, since the original chain $\{X_n \}$ is irreducible, the embedded chain $\{Y_n \}$ is also irreducible, since any finite-length positive-probability path from state $i\in F$ to state $j\in F$ in chain $\{X_n \}$ corresponds to a valid path between the two states in chain $\{Y_n \}$ when the states visited outside of set $F$ are removed.

    Define $\Tilde{\tau}_z \coloneqq \inf \{n \ge 0: Y_n  =z\}$. Noting that $Y_n$ is well-defined almost surely for all $n \in \mathbb{N}$, the following inequality holds pointwise in the extended real sense, almost surely:
    \begin{equation}
        \sum_{t=0}^{\tau_z-1} \mathds{I}\{X_t \in F \} \le \mathds{I}\{X_0 \in F\} + \Tilde{\tau}_z,\,\, \, \, \mathbb{P}_x-\mathrm{a.s}.
    \end{equation}
Taking expectation, $\mathbb{E}_x\left[\sum_{t=0}^{\tau_z-1} \mathds{I}\{X_t \in F \}\right] \le \mathds{I}\{x \in F\} + \mathbb{E}_x[\Tilde{\tau}_z]$. Hence, we now need to show that there exists a uniform bound for $\mathbb{E}_x[\Tilde{\tau}_z]$ across $\mathcal{X}$.

Noting that $Y_0=X_{\tau_F^{(1)}}$ is well-defined almost surely for any initial state $x\in \mathcal{X}$, we can exploit the strong Markov property to write $\mathbb{E}_x[\Tilde{\tau}_z]=\sum_{y\in F}\mathbb{E}_y^{Y}[\Tilde{\tau}_z]\,\mathbb{P}_x(Y_0=y)$, where $\mathbb{E}_y^Y[\cdot]$ corresponds to the expectation taken over the probability law of embedded chain starting from $Y_0=y$. It is well-established in Markov chain theory that, for irreducible and finite-state Markov chains, the expected hitting time to a certain state from a given initial state is finite for all possible initial states. Hence, $\mathbb{E}_y^Y[\Tilde{\tau}_z ] < \infty \, \, \forall \,y\in F$ and thus $\max_{y \in F} \mathbb{E}_y^Y[\Tilde{\tau}_z ] < \infty$.

Therefore, if we define $L \coloneqq 1 + \max_{y \in F} \mathbb{E}_y^Y[\Tilde{\tau}_z ]$, then $\mathbb{E}_x\left[\sum_{t=0}^{\tau_z-1} \mathds{I}\{X_t \in F \}\right] \le L \, \, \forall \, x\in \mathcal{X}$. This concludes the proof.
\end{proof}

\textbf{Theorem 6.5}. \emph{As a consequence of (\ref{drift sqrt}), and Lemma 6.2, and Lemma 6.4, the following are satisfied for all $x \in \mathcal{X}$:}
\begin{gather}
    (i) \,\mathbb{E}_x[\tau_z]\le \mathbb{E}_x\left[\sum_{t=0}^{\tau_z-1} W^{1/2}(X_t)\right ]\le \dfrac{W^{1/2}(x)+K^{1/2}L-1}{1-\lambda^{1/2}}, \nonumber \\ (ii) \, \dfrac{\mathbb{E}_x\left[ \sum_{t=0}^{\tau_z-1}W^{1/2}(X_t)\right]}{\mathbb{E}_x[\tau_z]} \le  \dfrac{W^{1/2}(x)+K^{1/2}-1}{1-\lambda^{1/2}}.
\end{gather}

\begin{proof} 

Onto the inequality presented in Lemma 6.2, setting $\tau = \tau_z$, utilizing the Monotone Convergence Theorem and invoking Lemma 6.4 recovers the first statement of the theorem.

Revisiting Lemma 6.2, setting $\tau = \tau_z$, noting that $\mathbb{E}_x\left[ \sum_{t=0}^{(N \land \tau_z)-1}\mathds{I}\{ X_t \in F\}\right] \le \mathbb{E}_x[N \land \tau_z]$, and invoking the Monotone Convergence Theorem, we get
\begin{equation}
    \mathbb{E}_x[\tau_z]\le \mathbb{E}_x\left[\sum_{t=0}^{\tau_z-1}W^{1/2}(X_t)\right] \le \dfrac{W^{1/2}(x)-1 +K^{1/2}\,\mathbb{E}_x\left[ \tau_z\right]}{1-\lambda^{1/2}} .
\end{equation}

Dividing both sides of the right-hand inequality by $\mathbb{E}_x[\tau_z]$, we get
\begin{equation}
    \dfrac{\mathbb{E}_x\left[ \sum_{t=0}^{\tau_z-1}W^{1/2}(X_t)\right]}{\mathbb{E}_x[\tau_z]} \le  \dfrac{W^{1/2}(x)-1}{(1-\lambda^{1/2})\mathbb{E}_x[\tau_z]} + \dfrac{K^{1/2}}{1-\lambda^{1/2}} \le \dfrac{W^{1/2}(x) +K^{1/2}-1}{1-\lambda^{1/2}},
\end{equation}

where the last inequality follows from $\mathbb{E}_x[\tau_z]\ge 1$. This completes the proof.\end{proof}

\textbf{Theorem 6.7.} \emph{Owing to Assumption 3.3, Theorem 6.5, Assumption 6.6 and \cite[Theorem 6.4]{makowski1994poisson}, there exists a valid solution to the Poisson equation for our problem as follows:}
\begin{align}
\label{sol to pois eq hitting time 2}
    V_\theta(x)&= \mathbb{E}_x\left[ \sum_{t=0}^{\tau_z-1}\bigl(g(\theta,X_t)-\bar{g}(\theta)\bigr)
     \right]\nonumber\\
    &=\mathbb{E}_x\left[ \sum_{t=0}^{\tau_z-1}g(\theta,X_t)
     \right]-\dfrac{\mathbb{E}_x[\tau_z]}{\mathbb{E}_z[\tau_z]}\,\mathbb{E}_z\left[ \sum_{t=0}^{\tau_z-1}g(\theta,X_t)
     \right]\, \,\,\forall \, x\in\mathcal{X}, \,\, \forall \,\theta\in \widetilde{\Theta}.
\end{align}

\emph{Moreover, statements $(i)$ and $(ii)$ of Assumption 3.4 hold with the functions $B_1(x)$, $B_2(x)$ selected as follows:}
\begin{gather}
    (i) \,B_1(x)\coloneqq D_1\dfrac{(W^{1/2}(z)+K^{1/2}-\lambda^{1/2})(W^{1/2}(x)+K^{1/2}L-1)}{(1-\lambda^{1/2})^2}, \nonumber\\
    (ii) \, B_2(x) \coloneqq D_2\dfrac{(W^{1/2}(z)+K^{1/2}-\lambda^{1/2})(W^{1/2}(x)+K^{1/2}L-1)}{(1-\lambda^{1/2})^2}.
\end{gather}
\begin{proof}
     Statement $(i)$ of Theorem 6.5 implies the recurrence condition of \cite{makowski1994poisson}. For the integrability condition, note that every element $i$ of the vector $g(\theta,x)$ satisfies $|g(\theta,x)_i| \le \lVert g(\theta,x)\rVert_{\infty}\le N \lVert g(\theta,x)\rVert_*\le NA_1(x)\,$ for some $N \in (0,\infty)$ due to norm equivalence. Hence, through Assumptions 3.3 and 6.6 and statement $(i)$ of Theorem 6.5, we get the following for all $\theta \in \widetilde{\Theta}, \, x\in \mathcal{X}$:
     \begin{equation}
         \mathbb{E}_x\left[\sum_{t=0}^{\tau_z -1} |g(\theta,X_t)_i|\right] 
         \le N \,\mathbb{E}_x\left[\sum_{t=0}^{\tau_z -1} A_1(X_t)\right] \le ND_1\,\mathbb{E}_x\left[\sum_{t=0}^{\tau_z -1} W^{1/2}(X_t)\right] < \infty.
     \end{equation}

    Therefore, given any $\theta\in\widetilde{\Theta}$ and index $i$, the integrability condition of \cite{makowski1994poisson} holds with $c(x)=g(\theta,x)_i$, for all $x \in \mathcal{X}$. Note that, through the verification of these two conditions, we also get that $\lVert\bar{g}(\theta)\rVert_* \le \mathbb{E}_{X\sim\pi}\left[ \,\lVert g(\theta,X) \rVert_* \right] < \infty$ for all $\theta \in \widetilde{\Theta}$. Applying \cite[Theorem 6.4]{makowski1994poisson} to every element of $V_{\theta}(x)$ belonging to (\ref{eq4}), we get the claimed expression for $V_\theta(x)$.

    Applying the norm operator, Assumptions 3.3 and 6.6, and Theorem 6.5 onto the expression for $V_\theta(x)$, we get the following:

\begin{align}
        \lVert V_\theta(x) \rVert_* &= \left\lVert \mathbb{E}_x\left[\sum_{t=0}^{\tau_z-1} g(\theta,X_t)\right]- \dfrac{\mathbb{E}_x[\tau_z]}{\mathbb{E}_z[\tau_z]}\,\mathbb{E}_z\left[ \sum_{t=0}^{\tau_z-1}g(\theta,X_t)
        \right]  \right\rVert_* \nonumber \\
        &\le \mathbb{E}_x\left[\sum_{t=0}^{\tau_z-1} \left\lVert g(\theta,X_t)\right\rVert_* \right]+ \dfrac{\mathbb{E}_x[\tau_z]}{\mathbb{E}_z[\tau_z]}\,\mathbb{E}_z\left[ \sum_{t=0}^{\tau_z-1}\lVert g(\theta,X_t) \rVert_* \right] \nonumber \\
        &\le \mathbb{E}_x\left[\sum_{t=0}^{\tau_z-1}  A_1(X_t)\right]+ \dfrac{\mathbb{E}_x[\tau_z]}{\mathbb{E}_z[\tau_z]}\,\mathbb{E}_z\left[ \sum_{t=0}^{\tau_z-1}A_1(X_t) \right] \nonumber \\
        &\le D_1 \left\{\mathbb{E}_x\left[\sum_{t=0}^{\tau_z-1}  W^{1/2}(X_t)\right]+ \dfrac{\mathbb{E}_x[\tau_z]}{\mathbb{E}_z[\tau_z]}\,\mathbb{E}_z\left[ \sum_{t=0}^{\tau_z-1}W^{1/2}(X_t) \right] \right\}
    \end{align}

\begin{align}
    V_{\theta'}(x)-V_\theta(x) &= \mathbb{E}_x\left[\sum_{t=0}^{\tau_z-1} [g(\theta', X_t)-g(\theta,X_t)]\right]- \dfrac{\mathbb{E}_x[\tau_z]}{\mathbb{E}_z[\tau_z]}\,\mathbb{E}_z\left[ \sum_{t=0}^{\tau_z-1}[g(\theta',X_t)-g(\theta,X_t)]
        \right] \nonumber \\
    &= \mathbb{E}_x\left[\sum_{t=0}^{\tau_z-1} [\tilde{g}(\theta', X_t)-\tilde{g}(\theta,X_t)]\right]- \dfrac{\mathbb{E}_x[\tau_z]}{\mathbb{E}_z[\tau_z]}\,\mathbb{E}_z\left[ \sum_{t=0}^{\tau_z-1}[\tilde{g}(\theta',X_t)-\tilde{g}(\theta,X_t)]
        \right] \nonumber \\
    &\,\,\,\,\,\,\,\,+ (r(\theta')-r(\theta))\left(\mathbb{E}_x[\tau_z]-\dfrac{\mathbb{E}_x[\tau_z]}{\mathbb{E}_z[\tau_z]}\mathbb{E}_z[\tau_z] \right) \nonumber \\
    &= \mathbb{E}_x\left[\sum_{t=0}^{\tau_z-1} [\tilde{g}(\theta', X_t)-\tilde{g}(\theta,X_t)]\right]- \dfrac{\mathbb{E}_x[\tau_z]}{\mathbb{E}_z[\tau_z]}\,\mathbb{E}_z\left[ \sum_{t=0}^{\tau_z-1}[\tilde{g}(\theta',X_t)-\tilde{g}(\theta,X_t)]
        \right] \nonumber 
\end{align}
\begin{align}
        \implies\lVert V_{\theta'}(x)-V_\theta(x) \rVert 
        &\le \mathbb{E}_x\left[\sum_{t=0}^{\tau_z-1} \left\lVert \tilde{g}(\theta', X_t)-\tilde{g}(\theta,X_t)\right\rVert\right]+ \dfrac{\mathbb{E}_x[\tau_z]}{\mathbb{E}_z[\tau_z]}\,\mathbb{E}_z\left[ \sum_{t=0}^{\tau_z-1}\lVert \tilde{g}(\theta',X_t)-\tilde{g}(\theta,X_t) \rVert \right] \nonumber \\
        &\le \lVert \theta' - \theta\rVert \left\{\mathbb{E}_x\left[\sum_{t=0}^{\tau_z-1}  A_2(X_t)\right]+ \dfrac{\mathbb{E}_x[\tau_z]}{\mathbb{E}_z[\tau_z]}\,\mathbb{E}_z\left[ \sum_{t=0}^{\tau_z-1}A_2(X_t) \right] \right\} \nonumber \\
        &\le D_2\, \lVert \theta' - \theta\rVert \left\{\mathbb{E}_x\left[\sum_{t=0}^{\tau_z-1}  W^{1/2}(X_t)\right]+ \dfrac{\mathbb{E}_x[\tau_z]}{\mathbb{E}_z[\tau_z]}\,\mathbb{E}_z\left[ \sum_{t=0}^{\tau_z-1}W^{1/2}(X_t) \right] \right\}
    \end{align}

Through Theorem 6.5, we bound the common factor in the above two expressions as follows:
\begin{align}
    \mathbb{E}_x\left[\sum_{t=0}^{\tau_z-1}  W^{1/2}(X_t)\right]+ \dfrac{\mathbb{E}_x[\tau_z]}{\mathbb{E}_z[\tau_z]}\,\mathbb{E}_z\left[ \sum_{t=0}^{\tau_z-1}W^{1/2}(X_t) \right] &\le \left(1 + \dfrac{\mathbb{E}_z\left[ \sum_{t=0}^{\tau_z-1}W^{\frac12}(X_t) \right]}{\mathbb{E}_z[\tau_z]} \right)\mathbb{E}_x\left[\sum_{t=0}^{\tau_z-1}  W^{1/2}(X_t)\right] \nonumber \\
    &\le \dfrac{(W^{1/2}(z)+K^{1/2}-\lambda^{1/2})(W^{1/2}(x)+K^{1/2}L-1)}{(1-\lambda^{1/2})^2}
\end{align}

Hence, with $B_i(x)\coloneqq D_i\dfrac{(W^{1/2}(z)+K^{1/2}-\lambda^{1/2})(W^{1/2}(x)+K^{1/2}L-1)}{(1-\lambda^{1/2})^2}$ for $i \in \{1,2 \}$, the following inequalities are satisfied for all $\theta, \theta'\in \widetilde{\Theta}$ and $x \in \mathcal{X}$:
\begin{gather}
    (i)\,\,\lVert V_{\theta}(x)\rVert_* \le B_1(x), \nonumber \\
    (ii)\,\,\lVert V_{\theta'}(x)-V_{\theta}(x) \rVert \le B_2(x)\lVert\theta' - \theta \rVert.
\end{gather}
\end{proof}

\textbf{Theorem 6.8.} \emph{Assumption 3.5 holds with $C_1(x), C_2(x), C_3(x), C_4(x)$ selected as follows:}
    \begin{gather}
        (i)\,C_1(x) \coloneqq D_1\dfrac{(W^{1/2}(z)+K^{1/2}-\lambda^{1/2})(W^{1/2}(x)+\frac{K^{1/2}}{1-\lambda^{1/2}}+K^{1/2}L-1)}{(1-\lambda^{1/2})^2} \nonumber \\
        (ii)\,C_2(x)\coloneqq D_1^2\left(W(x)+\frac{K}{1-\lambda}\right) \nonumber \\ 
        (iii)\, C_3(x) \coloneqq D_1^2\dfrac{(W^{1/2}(z)+K^{1/2}-\lambda^{1/2})(W(x)+\frac{K}{1-\lambda}+(K^{1/2}L-1)(W^{1/2}(x)+\frac{K^{1/2}}{1-\lambda^{1/2}}))}{(1-\lambda^{1/2})^2} \nonumber \\
        (iv) \,C_4(x) \coloneqq D_1D_2\dfrac{(W^{1/2}(z)+K^{1/2}-\lambda^{1/2})(W(x)+\frac{K}{1-\lambda}+(K^{1/2}L-1)(W^{1/2}(x)+\frac{K^{1/2}}{1-\lambda^{1/2}}))}{(1-\lambda^{1/2})^2} \nonumber\\
    \end{gather}

\begin{proof}
    The quantities of interest are bounded by expectations of functions in terms of $W(\cdot)$ and $W^{1/2}(\cdot)$, as illustrated below:
    \begin{align}
        (i)\,&\mathbb{E}_{x_0}[B_1(X_t)] = D_1\dfrac{(W^{1/2}(z)+K^{1/2}-\lambda^{1/2})(\mathbb{E}_{x_0}[W^{1/2}(X_t)]+K^{1/2}L-1)}{(1-\lambda^{1/2})^2} \nonumber \\
        (ii)\,&\mathbb{E}_{x_0}[A_1(X_t)^2]\le D_1^2 \,\mathbb{E}_{x_0}[W(X_t)] \nonumber\\
        (iii) \,&\mathbb{E}_{x_0}[A_1(X_t)B_1(X_{t+1})]\nonumber \\&\le D_1^2\dfrac{(W^{1/2}(z)+K^{1/2}-\lambda^{1/2})(\mathbb{E}_{x_0}[W^{1/2}(X_t)W^{1/2}(X_{t+1})]+(K^{1/2}L-1)\mathbb{E}_{x_0}[W^{1/2}(X_t)])}{(1-\lambda^{1/2})^2} \nonumber\\
        (iv) \,&\mathbb{E}_{x_0}[A_1(X_t)B_2(X_{t+1})] \nonumber\\&\le D_1D_2\dfrac{(W^{1/2}(z)+K^{1/2}-\lambda^{1/2})(\mathbb{E}_{x_0}[W^{1/2}(X_t)W^{1/2}(X_{t+1})]+(K^{1/2}L-1)\mathbb{E}_{x_0}[W^{1/2}(X_t)])}{(1-\lambda^{1/2})^2} \nonumber\\
    \end{align}
    
Let $\{X_t\}_{t\in\mathbb{N}}$ be a sample path of the Markov chain, with initial state $X_0=x_0$ for some $x_0\in\mathcal{X}$. Using the assumed drift inequalities, bounding $\mathbb{E}_{x_0}[W(X_t)]$ and $\mathbb{E}_{x_0}[W^{1/2}(X_t)]$ is straightforward. A one-step recursion on the conditional expectation of $W(X_t)$ yields
\begin{equation}
    \mathbb{E}_{x_0}[W(X_t)]=\mathbb{E}_{x_0}\bigl[\,\mathbb{E}_{x_0}[W(X_t)|\mathcal{F}_{t-1}]\,\bigr] \le \lambda\,\mathbb{E}_{x_0}\left[ W(X_{t-1}) \right] + K
\end{equation}

Applying the recursion down to $X_0$, we get
\begin{equation}
    \mathbb{E}_{x_0}[W(X_t)] \le \lambda^t\,W(x_0)+K\frac{1-\lambda^t}{1-\lambda} \le W(x_0)+\frac{K}{1-\lambda}
\end{equation}

Replacing $W(\cdot), \lambda, K$ by $W^{1/2}(\cdot), \lambda^{1/2},$ and $K^{1/2}$, respectively, we get the corresponding relation for $W^{1/2}(\cdot):$
\begin{equation}
    \mathbb{E}_{x_0}[W^{1/2}(X_t)] \le W^{1/2}(x_0)+\frac{K^{1/2}}{1-\lambda^{1/2}}
\end{equation}

Bounding $\mathbb{E}_{x_0}[W^{1/2}(X_t)W^{1/2}(X_{t+1})]$ is straightforward, using the Cauchy-Schwarz inequality for nonnegative random variables:
\begin{equation}
    \mathbb{E}_{x_0}[W^{1/2}(X_t)W^{1/2}(X_{t+1})]\le \left( \mathbb{E}_{x_0}[W(X_t)]\,\mathbb{E}_{x_0}[W(X_{t+1})]\right)^{1/2} \le W(x_0)+\dfrac{K}{1-\lambda}.
\end{equation}

Hence, with $C_1(x),C_2(x),C_3(x),C_4(x)$ defined as in the statement of the theorem, the following inequalities hold for all $t \in \mathbb{N}$, $x_0 \in \mathcal{X}$:
\begin{gather}
    (i)\,\,\mathbb{E}_{x_0}[B_1(X_t)]\le C_1(x_0),\nonumber\\
    (ii)\,\,\mathbb{E}_{x_0}[A_1(X_t)^2] \le C_2(x_0), \nonumber \\
    (iii)\,\,\mathbb{E}_{x_0}[A_1(X_t)B_1(X_{t+1})] \le C_3(x_0), \nonumber \\
    (iv)\,\,\mathbb{E}_{x_0}[A_1(X_t)B_2(X_{t+1})] \le C_4(x_0). 
\end{gather}
\end{proof}

\end{document}